\documentclass[10pt,twoside,a4paper]{article}

\usepackage[latin1]{inputenc}
\usepackage[T1]{fontenc}
\usepackage{amsmath,amssymb,amsthm}
\usepackage{array}
\usepackage{color}
\usepackage{graphicx}

\newcommand{\Id}{\mathbf{1}}

\renewcommand{\>}{\rangle}
\newcommand{\<}{\langle}

\newcommand{\IM}{\operatorname{\mathcal{M}}}
\newcommand{\BB}{\operatorname{\mathcal{B}}}
\newcommand{\EE}{\operatorname{\mathcal{E}}}
\newcommand{\DM}{\overline{\operatorname{\mathcal{M}}}}
\newcommand{\IA}{\operatorname{\mathfrak{A}}}

\renewcommand{\SS}{\operatorname{\mathfrak{S}}}

\newcommand{\Si}{\dot{S}}

\renewcommand{\H}{\mathcal{H}}
\newcommand{\J}{\mathcal{J}}
\newcommand{\cyl}{\operatorname{cyl}}
\newcommand{\CP}{\mathbb{C}\mathbb{P}^1}

\newcommand{\CR}{\overline{\partial}}
\newcommand{\uz}{\underline{z}}
\newcommand{\Mph}{M_{\phi}}

\newcommand{\Aut}{\operatorname{Aut}}
\newcommand{\Symp}{\operatorname{Symp}}

\newcommand{\Sp}{\operatorname{Sp}}

\newcommand{\Crit}{\operatorname{Crit}}

\newcommand{\cst}{\operatorname{const}}
\newcommand{\st}{\operatorname{st}}

\newcommand{\loc}{\operatorname{loc}}
\newcommand{\reg}{\operatorname{reg}}
\newcommand{\End}{\operatorname{End}}

\newcommand{\hb}{\bar{h}}

\newcommand{\diag}{\operatorname{diag}}
\newcommand{\ind}{\operatorname{ind}}
\newcommand{\coker}{\operatorname{coker}}
\newcommand{\del}{\partial}

\newcommand{\RS}{\IR \times S^1}
\newcommand{\Ju}{\underline{J}}

\renewcommand{\qed}{\square}

\newcommand{\IC}{\operatorname{\mathbb{C}}}
\newcommand{\IZ}{\operatorname{\mathbb{Z}}}
\newcommand{\IQ}{\operatorname{\mathbb{Q}}}
\newcommand{\IR}{\operatorname{\mathbb{R}}}
\newcommand{\IN}{\operatorname{\mathbb{N}}}
\newcommand{\kom}[1]{}

\newtheorem{maintheorem}{Main Theorem}
\newtheorem{theorem}{Theorem}[section]
\newtheorem{proposition}[theorem]{Proposition}
\newtheorem{definition}[theorem]{Definition}
\newtheorem{lemma}[theorem]{Lemma}
\newtheorem{corollary}[theorem]{Corollary}

\author{Oliver Fabert\\Mathematisches Institut der Ludwig-Maximilians-Universit\"at M\"unchen
        \\Theresienstr. 39, 80333 M\"unchen, Germany\\fabert@math.lmu.de}
\title{Contact Homology of Hamiltonian mapping tori}
\begin{document}

\maketitle

\begin{abstract}
In the general geometric setup for symplectic field theory the contact manifolds 
can be replaced by mapping tori $\Mph$ of symplectic manifolds $(M,\omega)$ with symplectomorphisms $\phi$. 
While the cylindrical contact homology of $\Mph$ is given by the Floer homologies of powers of $\phi$, 
the other algebraic invariants of symplectic field theory for $\Mph$ provide natural generalizations of symplectic Floer homology.   
For symplectically aspherical $M$ and Hamiltonian $\phi$ we study the moduli spaces of rational curves and prove a transversality result, 
which does {\it not} need the polyfold theory by Hofer, Wysocki and Zehnder. We use our result to compute the full contact homology of $M_{\phi}\cong S^1\times M$. 
\end{abstract}

\tableofcontents

\section{Introduction and main results}

\subsection{Symplectic field theory in the Floer case}
This paper is concerned with symplectic field theory in the Floer case. Symplectic field theory (SFT) is a very 
large project designed to describe in a unified way the theory of holomorphic curves in symplectic and contact topology. To be more precise, it approaches Gromov-Witten theory in the spirit of a topological quantum field theory by counting holomorphic curves in cylinders over contact manifolds and symplectic cobordisms between them. It was initiated by Eliashberg, Givental and Hofer in their paper [EGH] and since then has found many striking applications in symplectic geometry and beyond. While most of the current applications lie in finding invariants for contact manifolds, there exists a generalized geometric setup for symplectic field theory, 
which contains contact manifolds as special case. \\
    
Following [BEHWZ] and [CM2] a {\it Hamiltonian structure} on a closed $(2m-1)$-dimensional manifold $V$ is a closed two-form $\omega$ 
on $V$ which is maximally nondegenerate in the sense that $\ker\omega=\{v\in TV:\omega(v,\cdot)=0\}$ 
is a one-dimensional distribution. Note that here we (and [CM2]) differ slightly from [EKP]. 
The Hamiltonian structure is required to be {\it stable} in the sense that there exists a 
one-form $\lambda$ on $V$ such that $\ker\omega\subset\ker d\lambda$ and $\lambda(v)\neq 0$ for all $v\in\ker\omega-\{0\}$. 
Any stable Hamiltonian structure $(\omega,\lambda)$ defines a symplectic hyperplane distribution $(\xi=\ker\lambda,\omega_{\xi})$, 
where $\omega_{\xi}$ is the restriction of $\omega$, and 
a vector field $R$ on $V$ by requiring $R\in\ker\omega$ and $\lambda(R)=1$ which is called the Reeb vector field of the 
stable Hamiltonian structure. \\

Examples for closed manifolds $V$ with a stable Hamiltonian structure $(\omega,\lambda)$ 
are contact manifolds, circle bundles and mapping tori ([BEHWZ],[CM2]). 
For this note that when $\lambda$ is a contact form on $V$, then it is easy to check that $(\omega:=d\lambda,\lambda)$ is a stable 
Hamiltonian structure and the symplectic hyperplane distribution agrees with the contact structure. 
For the other two cases, let $(M,\omega)$ be a symplectic manifold. Then any principal circle 
bundle $S^1\to V\to M$ and any symplectic mapping torus $M\to V\to S^1$, i.e., $V=\Mph=\IR\times M/\{(t,p)\sim(t+1,\phi(p))\}$ for 
$\phi\in\Symp(M,\omega)$ carries also a stable Hamiltonian structure. For the circle bundle the Hamiltonian 
structure is given by the pullback $\pi^*\omega$ under the bundle projection and the one-form $\lambda$ is given by any $S^1$-connection form. On the other hand, the stable Hamiltonian structure on the mapping torus $V=\Mph$ is given by lifting the symplectic form to $\omega\in\Omega^2(\Mph)$ via the natural flat connection $TV=TS^1\oplus TM$ and setting $\lambda=dt$ for the natural $S^1$-coordinate $t$ on $\Mph$. While in the mapping torus case $\xi$ is always integrable, in the circle bundle case the hyperplane distribution $\xi$ may be integrable or non-integrable, even contact. \\ 

Symplectic field theory assigns algebraic invariants to closed manifolds $V$ with a stable Hamiltonian structure. 
The invariants are defined by counting $\Ju$-holomorphic curves in $\IR\times V$ with finite energy, 
where the underlying closed Riemann surfaces are explicitly allowed to have punctures, i.e., single points are removed. 
The almost complex structure $\Ju$ on the cylindrical manifold $\IR\times V$ is required to be {\it cylindrical} in the sense that it is  
$\IR$-independent, links the two natural vector fields on $\IR\times V$, namely the Reeb vector field 
$R$ and the $\IR$-direction $\del_s$, by $\Ju\del_s=R$, and turns the symplectic hyperplane distribution on $V$ into a complex subbundle of $TV$, $\xi=TV\cap \Ju TV$. It follows that a cylindrical almost complex structure $\Ju$ on $\IR\times V$ is determined by its restriction 
$\Ju_{\xi}$ to $\xi\subset TV$, which is required to be $\omega_{\xi}$-compatible in the sense that 
$\omega_{\xi}(\cdot,\Ju_{\xi}\cdot)$ defines a metric on $\xi$. Note that in [CM2] such almost complex structures $\Ju$ are called compatible 
with the stable Hamiltonian structure and that the set of these almost complex structures is non-empty and contractible. \\

While the punctured curves in symplectic field theory may have arbitrary genus and arbitrary numbers of positive and negative punctures, it  
is shown in [EGH] that there exist algebraic invariants counting only special types of curves: While in {\it rational symplectic field theory} one counts punctured curves with genus zero, {\it contact homology} is defined by further restricting to punctured spheres with only one positive puncture. Further restricting to spheres with both just one negative and one positive puncture, i.e., cylinders, the resulting algebraic invariant is called {\it cylindrical contact homology}. Note however that contact homology and cylindrical contact homology are not always defined. In order to prove the well-definedness of (cylindrical) contact homology it however suffices to show that there are no punctured holomorphic curves where all punctures are negative (or all punctures are positive). While the existence of holomorphic curves without positive punctures can be excluded for all contact manifolds using the maximum principle, which shows that contact homology is well-defined for all contact manifolds, it can be seen from homological reasons that for mapping tori $M_{\phi}$ there cannot exist holomorphic curves in $\IR\times M_{\phi}$ carrying just one type of punctures, which shows that in this case both contact homology and cylindrical contact homology are defined. \\

Symplectic field theory hence provides a wealth of invariants. However, almost all computations performed so far only use the simplest one, cylindrical contact homology: While cylindrical contact homology is computed e.g. for subcritical Stein-fillable contact manifolds ([Y1]), Brieskorn varieties ([K]) and toroidal three-manifolds ([BC]), computations of the higher invariants are performed so far only for overtwisted contact manifolds in [Y2] and sketched in [EGH] for prequantization spaces and in [CL] for unit cotangent bundle of tori.   

\subsection{Main theorem and outline of the proof}
While it can be seen that the cylindrical contact homology for mapping tori $\Mph$ agrees with the Floer homology of the powers of $\phi$, i.e., the subcomplex for the period $T\in\IN$ agrees with the Floer homology of $\phi^T$, the other algebraic invariants of symplectic field theory, in particular, the full contact homology, provide natural generalizations of symplectic Floer homology. While Floer homology for Hamiltonian symplectomorphisms over a suitable coefficient ring is known to be isomorphic to the tensor product of the singular homology with rational coefficients of the underlying symplectic manifold with the graded group algebra $\IQ[H_2(M)]$ generated by $H_2(M)$,
\begin{equation*} \IQ[H_2(M)] = \{\sum q(A) e^A: A\in H_2(M), q(A) \in \IQ\}, \;\;\deg e^A = \<c_1(TM),A\>, \end{equation*} 
the case of arbitrary symplectomorphisms is much more complicated, see [CC] and the references therein. So we restrict our attention to the Hamiltonian case, where the symplectomorphism $\phi$ is Hamiltonian, i.e., the time-one map of the symplectic flow of a Hamiltonian $H: S^1\times M \to \IR$. In this case the Hamiltonian flow $\phi^H$ provides us with a natural diffeomorphism $\Mph\cong S^1 \times M$, so that we can replace $\Mph$ by $S^1\times M$ equipped with the pullback stable Hamiltonian structure $(\omega^H,\lambda^H)$ on $S^1\times M$ given by $\omega^H = \omega + dH\wedge dt$, $\lambda^H=dt$ with symplectic bundle $\xi^H=TM$ and Reeb vector field $R^H = \del_t + X^H_t$, where $X^H_t$ is the symplectic gradient of $H_t=H(t,\cdot)$. In [EKP] this 
is also called the {\it Floer case.} Furthermore $(\IR\times\Mph,\Ju)$ can be identified with $(\RS\times M,\Ju^H)$ equipped with the pullback cylindrical almost 
complex structure, which is nonstandard in the sense that the splitting $T(\RS\times M) = \IR^2\oplus TM$ is not $\Ju^H$-complex. Let $HC_*(\Mph,\Ju)$ denote the contact homology of the symplectic mapping torus $\Mph$ with chosen cylindrical almost complex structure $\Ju$ on $\IR\times\Mph$. \\

\begin{maintheorem} Let $(M,\omega)$ be a closed symplectic manifold, which is symplectically aspherical, $\<[\omega],\pi_2(M)\>=0$, and let $\phi:M\to M$ be a Hamiltonian symplectomorphism. Then we have 
\begin{eqnarray*}
 HC_*(\Mph,\Ju) &\cong& \SS\bigl(\bigoplus_{\IN} H_{*-2}(M,\IQ)\bigr) \otimes \IQ[H_2(M)], 
\end{eqnarray*}
where $\SS$ is the graded symmetric algebra functor. \end{maintheorem} 

For the proof we observe that the cylindrical almost complex structure $\Ju^H$ on $\RS\times M$ is specified by the choice of an $S^1$-family of 
almost complex structures $J_t$ on $M$ and an $S^1$-dependent Hamiltonian $H: S^1\times M\to \IR$. 
In order to get an $S^1$-symmetry on moduli spaces of curves with three or more punctures, we restrict ourselves 
to almost complex structures $J_t$ and Hamiltonians $H_t$, which are independent of $t\in S^1$, so that only holomorphic cylinders 
need to be counted for the differential in contact homology. \\
  
We achieve transversality for all moduli spaces by considering {\it domain-dependent Hamiltonian perturbations}. 
This means that, for defining the Cauchy-Riemann operator for curves, we allow the Hamiltonian to depend explicitly 
on points on the punctured sphere underlying the curve whenever 
the punctured sphere is stable, i.e., there are no nontrivial automorphisms, where we follow the ideas in [CM1]. Note however that 
in contrast to the Gromov-Witten case we now have to make coherent choices for the different moduli spaces 
{\it simultaneously}, i.e., the different Hamiltonian perturbations must be compatible with gluing of curves in 
symplectic field theory. For the cylindrical moduli spaces the Hamiltonian perturbation is domain-independent, and it 
is known from Floer theory that in general we must allow $H$ to depend 
explicitly on $t\in S^1$ to achieve nondegeneracy of the periodic orbits and transversality for the moduli spaces 
of Floer trajectories. However, the gluing compatibility requires that also 
the Hamiltonian perturbation for the cylindrical moduli spaces is independent of $t\in S^1$. 
We solve this problem by considering Hamiltonians $H$, which are so small in the $C^2$-norm that all orbits are critical points of $H$ and all cylinders between these orbits correspond to gradient flow lines between the underlying critical points. Note however that we cannot achieve this with a single Hamiltonian function, but have to rescale the function depending on the period $T\in\IN$, which in turn implies that we have to compute the contact homology using an infinite sequence of different Hamiltonian functions. \\

Observe that the closed orbits of the Reeb vector field $R^H$ on $S^1\times M$ have integer periods, where the set of closed orbits of period $T\in\IN$ is naturally identified with the $T$-periodic orbits of $X^H$ on $M$. It follows that the chain complex $(\IA,\del)$ for contact homology naturally splits, $\IA=\bigoplus_{T\in\IN}\IA^T$, where $\IA^T$ is generated by all monomials $q_{(x_1,T_1)}...q_{(x_n,T_n)}$, with $T_i$-periodic orbits $(x_i,T_i)$ and $T_1+...+T_n=T$, and it is easily seen from homological reasons that this splitting is respected by the differential $\del$. Furthermore, given two different Hamiltonian functions $H_1,H_2: S^1\times M\to\IR$ the corresponding chain map $\Phi:(\IA_1,\del_1)\to (\IA_2,\del_2)$, defined as in [EGH] by counting holomorphic curves in $\RS\times M$ equipped with a non-cylindrical almost complex structure $\Ju^{\tilde{H}}$, which itself can be defined using a homotopy $\tilde{H}:\RS\times M\to\IR$ from $H_1$ to $H_2$, also respects the splittings $\IA_1=\bigoplus_{T\in\IN} \IA^T_1$, $\IA_2=\bigoplus_{T\in\IN} \IA^T_2$. \\

Let $T_N\in\IN$ be a sequence of (maximal) periods with $T_N\leq T_{N+1}$ and $\lim_{N\to\infty} T_N=\infty$ and 
let $H_N:S^1\times M\to\IR$, $N\in\IN$ be a sequence of Hamiltonians with corresponding chain complexes $(\IA_N,\del_N)$, $N\in\IN$. Assume that for every $N\in\IN$ we have defined a chain map $\Phi_N: (\IA_N,\del_N)\to (\IA_{N+1},\del_{N+1})$ using a homotopy $\tilde{H}_N:\RS\times M\to\IR$ interpolating between $H_N$ and $H_{N+1}$, which by the above arguments restricts to a map from $\IA^T_N$ to 
$\IA^T_{N+1}$ for every $T\in\IN$. Defining
\begin{equation*} 
 HC_*^{\leq T_N}(S^1\times M,\Ju^{H_N}) = H_*(\IA^{\leq T_N}_N,\del_N) = \bigoplus_{T\leq T_N} H_*(\IA^T_N,\del_N)
\end{equation*}
we obtain a directed system $(C_N,\Phi_{N,M})$ with $C_N=HC_*^{\leq T_N}(S^1\times M,\Ju^{H_N})$ and $\Phi_{N,M}=\Phi_N\circ\Phi_{N+1}\circ ...\circ\Phi_{M-1}\circ\Phi_M$ for $N\leq M$. Setting $T_N=2^N$ we prove the main result by showing that for every $S^1$-independent Hamiltonian $H:M\to\IR$, which is sufficiently small in the $C^2$-norm and Morse, there  
is an isomorphism 
\begin{eqnarray*}
 \lim_{N\to\infty} HC_*^{\leq 2^N}(S^1\times M,\Ju^{H/2^N}) &\cong& \SS\bigl(\bigoplus_{\IN} H_{*-2}(M,\IQ)\bigr) \otimes \IQ[H_2(M)]. 
\end{eqnarray*}

This paper is organized as follows. \\

While we prove in 2.1 all the fundamental results about pseudoholomorphic curves in Hamiltonian mapping tori, 
subsection 2.2 is devoted to explaining the central ideas of the proof the main theorem, namely how we get an $S^1$-symmetry on all moduli spaces of domain-stable curves, but still have nondegeneracy for the closed orbits and transversality for all moduli spaces. We collect all the important results about the moduli spaces in theorem 2.6. After recalling the definition of the Deligne-Mumford space of stable punctured spheres in 3.1, we define the underlying domain-dependent Hamiltonian perturbations in 3.2 and prove in 3.3 that the construction is compatible with the SFT compactness theorem. After describing in detail the necessary Banach manifold setup for our Fredholm problems in 4.1, 
we prove in 4.2 the fundamental transversality result for the Cauchy-Riemann operator. 
Since all our results only hold up to a maximal period for the asymptotic orbits, i.e., we 
have to rescale our Hamiltonian perturbation during the computation of contact homology in section 6, we generalize all 
our previous results to homotopies of Hamiltonian perturbations in 5.1 and 5.2. After describing the chain complex underlying contact homology in 6.1, 
we prove the main theorem using our previous results about moduli spaces of holomorphic curves in $\RS\times M$. \\

{\bf Acknowledgements} This research was supported by the German Research Foundation (DFG). 
The author thanks U. Frauenfelder, M. Hutchings and K. Mohnke for useful 
conversations and their interest in his work. Special thanks finally go to my advisor Kai Cieliebak and to 
Dietmar Salamon, who gave me the chance to stay at ETH Zurich for the winter term 2006/07, for their support. Finally thanks go to the referee 
for his valuable comments.  

\section{Moduli spaces}

\subsection{Holomorphic curves in $\RS\times M$}

Let $(M,\omega)$ be a closed symplectic manifold and let $\phi$ be a symplectomorphism on it. As already explained in the 
introduction, the corresponding mapping torus $\Mph=\IR\times M/\{(t,p)\sim(t+1,\phi(p))\}$ carries a natural stable 
Hamiltonian structure $(\omega,\lambda)$ given by lifting the symplectic form $\omega$ to a two-form on $\Mph$ via 
the flat connection $T\Mph=TS^1\oplus TM$ and setting $\lambda=dt$. It follows that the corresponding symplectic vector bundle 
$\xi=\ker\lambda$ is given by $TM$ and the Reeb vector field $R$ agrees with the $S^1$-direction $\del_t$ on $\Mph$. 
In this paper we restrict ourselves to the case where $\<[\omega],\pi_2(M)\> = 0$ and $\phi$ is Hamiltonian, i.e., 
the time-one map of the flow of a Hamiltonian $H: S^1\times M \to \IR$. In this case observe that the Hamiltonian flow 
$\phi^H$ provides us with the natural diffeomorphism 
\begin{equation*}
 \Phi: S^1 \times M \stackrel{\cong}{\longrightarrow} \Mph,\,(t,p) \mapsto (t,\phi^H(t,p)),  
\end{equation*}
so that we can replace $\Mph$ by $S^1\times M$ equipped with the pullback stable Hamiltonian structure. \\
\\
\begin{proposition} The pullback stable Hamiltonian structure $(\omega^H,\lambda^H)$ on $S^1\times M$ is given by 
\begin{equation*} \omega^H \;=\; \omega + dH\wedge dt,\;\;\; \lambda^H\;=\;dt \end{equation*} 
with symplectic bundle $\xi^H$ and Reeb vector field $R^H$ given by 
\begin{equation*} \xi^H \;=\; TM,\;\;\; R^H\;=\;\del_t + X^H_t, \end{equation*} 
where $X^H_t$ is the symplectic gradient of $H_t=H(t,\cdot)$. \end{proposition}  
\noindent
{\it Proof:} Using 
\begin{equation*} d\Phi = (\Id,X^H_t\otimes dt + d\phi^H_t):  TS^1\oplus TM \to TS^1\oplus TM \end{equation*} 
we compute for $v_1=(v_{11},v_{12}),v_2=(v_{21},v_{22})\in TS^1\oplus TM$, 
\begin{eqnarray*} 
\omega^H(v_1,v_2) &=& \omega(d\Phi(v_1), d\Phi(v_2)) \\
                  &=& \omega((X^H_t\otimes dt)(v_{11}) + d\phi^H_t(v_{12}), (X^H_t\otimes dt)(v_{21}) + d\phi^H_t(v_{22}))\\ 
                  &=& \omega(X^H_t,X^H_t) dt(v_{11}) dt(v_{21}) + \omega(d\phi^H_t(v_{12}),d\phi^H_t(v_{22})) \\
                  && + \omega(X^H_t,d\phi^H_t(v_{22})) dt(v_{11}) + \omega(d\phi^H_t(v_{12}),X^H_t) dt(v_{21}) \\ 
                  &=& \omega(v_{12},v_{22}) + \omega(d\phi^H_t(v_{12}),X^H_t) dt(v_{21}) - \omega(d\phi^H_t(v_{22}),X^H_t) dt(v_{11}) \\
                  &=& \omega(v_1,v_2) + (dH\wedge dt)(v_1,v_2) 
\end{eqnarray*}
and $\lambda^H = \lambda \circ d\Phi = dt$. On the other hand, it directly follows that $\xi^H=TM$, while $R^H=\del_t-X^H_t$ spans 
the kernel of $\omega^H$, 
\begin{eqnarray*} 
\omega^H(\cdot,R^H) &=& \omega(\cdot,\del_t-X^H_t) + dH \cdot dt(\del_t+X^H_t) - dH(\del_t+X^H_t)\cdot dt \\
                    &=& - \omega(\cdot,X^H_t) + dH = 0 
\end{eqnarray*} 
with $\lambda^H(R^H) = dt(\del_t-X^H_t) = 1$. $\qed$ \\    

As in the introduction we consider an almost complex structure $\Ju$ on the cylindrical manifold $\RS\times M$, which 
is required to be cylindrical in the sense that it is $\IR$-independent, links the Reeb vector field $R^H$ and the $\IR$-direction $\del_s$, 
by $\Ju\del_s=R^H=\del_t+X^H_t$ and turns the symplectic hyperplane distribution $\xi^H=TM$ into a complex subbundle of $T(S^1\times M)$. 
It follows that $\Ju$ on $\IR\times S^1\times M$ is determined by its restriction to $\xi^H=TM$, which is required to be $\omega_{\xi^H}$-compatible, so 
that $\Ju$ is determined by the $S^1$-dependent Hamiltonian $H_t$ and an $S^1$-family of $\omega$-compatible almost complex structures 
$J_t$ on the symplectic manifold $(M,\omega)$. \\
 
Let us recall the definition of moduli spaces of holomorphic curves studied in rational SFT in the general setup. 
Let $(V,\omega,\lambda)$ be a closed manifold with stable Hamiltonian structure with symplectic hyperplane distribution $\xi$ and Reeb vector field $R$ and 
let $\Ju$ be a compatible cylindrical almost complex structure on $\IR\times V$. Let $P^+,P^-$ be two ordered sets of closed 
orbits $\gamma$ of the Reeb vector field 
$R$ on $V$, i.e., $\gamma: \IR\to V$, $\gamma(t+T)=\gamma(t)$, $\dot{\gamma}=R$, where $T>0$ denotes the period of $\gamma$. 
Then the (parametrized) moduli space $\IM^0(V;P^+,P^-,\Ju)$ consists of tuples $(F,(z_k^{\pm}))$, where 
$\{z^{\pm}_1,...,z^{\pm}_{n^{\pm}}\}$ are two disjoint ordered sets of points on $\CP$, 
which are called positive and negative punctures, respectively. The map $F: \Si \to \IR \times V$  starting 
from the punctured Riemann surface $\Si = \CP -  \{(z_k^{\pm})\}$ is required to satisfy the 
Cauchy-Riemann equation 
\begin{equation*}
 \CR_{\Ju} F = dF + \Ju(F) \cdot dF \cdot i = 0
\end{equation*}
with the complex structure $i$ on $\CP$. Assuming we have chosen cylindrical coordinates $\psi^{\pm}_k: \IR^{\pm} 
\times S^1 \to \Si$ around each puncture $z^{\pm}_k$ in the sense that $\psi_k^{\pm}(\pm\infty,t)=z_k^{\pm}$, 
the map $F$ is additionally required to show for all $k=1,...,n^{\pm}$ the 
asymptotic behaviour
\begin{equation*}
 \lim_{s\to\pm\infty} (F \circ \psi^{\pm}_k) (s,t+t_0) = 
 (\pm \infty,\gamma^{\pm}_k(T^{\pm}_kt))
\end{equation*}
with some $t_0\in S^1$ and the orbits $\gamma^{\pm}_k\in P^{\pm}$, where $T^{\pm}_k>0$ denotes period of $\gamma^{\pm}_k$.
Observe that the group $\Aut(\CP)$ of Moebius transformations 
acts on elements in $\IM^0(V;P^+,P^-,\Ju)$ in an obvious way, 
\begin{equation*} \varphi.(F,(z^{\pm}_k)) = (F\circ\varphi^{-1},\varphi(z^{\pm}_k)),\;\;\;\varphi\in\Aut(\CP), \end{equation*}
and we obtain the moduli space $\IM(V;P^+,P^-,\Ju)$ studied in symplectic field theory by dividing out this action. \\

It remains to identify the occuring objects in our special case. First, one immediately verifies that all closed orbits $\gamma$ of the 
vector field $R^H=\del_t-X^H_t$ on $S^1\times M$ are of the form 
\begin{equation*} \gamma(t)=(t+t_0,x(t)), \end{equation*} 
and therefore have natural numbers $T\in\IN$, i.e., the winding number around the $S^1$-factor, as periods. Since we study closed Reeb orbits up to 
reparametrization, we can set $t_0=0$, so that $\gamma$ can be identified with $x:\IR/T\IZ\to M$, which is a 
$T$-periodic orbit of the Hamiltonian vector field, 
\begin{equation*} \dot{x}(t)=X^H_t(x(t)). \end{equation*} 
Hence we will in the following write $\gamma=(x,T)$, where $T\in\IN$ is the period and $x$ is a $T$-periodic orbit of the Hamiltonian 
$H$. We denote the set of $T$-periodic orbits of the Reeb vector field $R^H$ on $S^1\times M$ by $P(H,T)$.
 
For the moduli spaces of curves observe that in $\RS\times M$ we can naturally write the holomorphic map $F$ as a product, 
\begin{equation*}
 F = (h,u): \Si \to (\RS) \times M\;.
\end{equation*} 
\\
\begin{proposition} $F:\Si\to\RS\times M$ is $\Ju$-holomorphic precisely when $h=(h_1,h_2):\Si\to\RS$ is holomorphic and 
$u:\Si\to M$ satisfies the $h$-dependent perturbed Cauchy-Riemann equation of Floer type, 
\begin{eqnarray*} 
\CR_{J,H,h} u &=& \Lambda^{0,1}(du + X^H(h_2,u)\otimes dh_2) \\ 
              &=& du + X^H(h_2,u)\otimes dh_2 + J(h_2,u)\cdot(du+X^H(h_2,u)\otimes dh_2) \cdot i.
\end{eqnarray*}
\end{proposition}

\noindent{\it Proof:} Observing that $\Ju(t,p): T(\RS)\oplus TM\to T(\RS)\oplus TM$ is given by 
\begin{equation*} \Ju(t,p)=\binom{\;\;i\;\;\;\;0\;\;}{\Delta(t,p)\;\;J_t(p)} \end{equation*} 
with $\Delta(t,p) = -X^H_t(p)\otimes ds + J_t(p) X^H_t(p)\otimes dt$ we compute 
\begin{eqnarray*} 
 && (dh,du) + \Ju(h,u) \cdot (dh,du) \cdot i \\
 &=& (dh + i \cdot dh \cdot i, \\
 && \;du + (J(h_2,u)\cdot du - X^H(h_2,u)\otimes dh_1 + J(h_2,u)X^H(h_2,u)\otimes dh_2) \cdot i) \\
 &=& (\CR h, du - X^H(h_2,u)\otimes dh_1\cdot i + J(h_2,u)\cdot(du + X^H(h_2,u)\otimes dh_2) \cdot i). 
\end{eqnarray*} 
Finally observe that $dh_1 \cdot i = -dh_2$ if $\CR h =0$. $\qed$ \\

Recalling that our orbit sets are given by $P^{\pm} = \{(x^{\pm}_1,T^{\pm}_1),...,(x^{\pm}_{n^{\pm}},T^{\pm}_{n^{\pm}})\}$, 
we use the rigidity of holomorphic maps to prove the following statement about the map component $h: \Si\to\RS$. Let 
$T^{\pm} = T^{\pm}_1+...+T^{\pm}_{n^{\pm}}$ denote the total period above and below, respectively. \\
\\
\begin{lemma} The map $h = (h_1,h_2)$ exists if and only if $T^+=T^-$ and is unique up a shift $(s_0,t_0) \in \RS$,
\begin{equation*}
   h(z) = h^0(z) + (s_0,t_0)
\end{equation*}
for some fixed map $h^0=(h^0_1,h^0_2)$. In particular, every holomorphic cylinder has a positive and a negative puncture, there are no holomorphic planes and all holomorphic spheres are constant. \end{lemma}
\noindent  
{\it Proof:} The asymptotic behavior of the map $F$ near the punctures implies that 
\begin{equation*}
h\circ\psi_k(s,t+t_0)\stackrel{s\to\pm\infty}{\longrightarrow}(\pm\infty,T_kt) 
\end{equation*} 
with some $t_0\in S^1$. Identifying $\RS \cong \CP - \{0,\infty\}$, it follows that 
$h$ extends to a meromorphic function $h$ on $\CP$ with $z^+_1,...,z^+_{n^+}$ poles of order $T^+_1,...,T^+_{n^+}$ and 
$z^-_{1},...,z^-_{n^-}$ zeros of order $T^-_1,...,T^-_{n^-}$. Since the zeroth Picard group of $\CP$ is trivial, i.e., every 
divisor of degree zero is a principal divisor, we get that such meromorphic functions exist precisely when $T^+=T^-$. On the 
other hand it follows from Liouville's theorem that they are uniquely determined up to a 
nonzero multiplicative factor, i.e., $h = a \cdot h^0$ with $a \in \IC^*\cong\RS$ for some fixed $h_0: \CP \to \CP$. 
For every $\Ju^H$-holomorphic sphere $(h,u)$ observe that $h$ is constant, $h=(s_0,t_0)$, and therefore $u$ is a $J_{t_0}$-holomorphic 
sphere in $M$, which must be constant by $\<[\omega],\pi_2(M)\>=0$. $\qed$ \\

Note that the lemma also holds when $\phi$ is no longer Hamiltonian by defining $h = \pi\circ F$ using the holomorphic bundle projection $\pi:\IR\times\Mph\to\RS$. \\

It follows that we only have to study punctured $\Ju^H$-holomorphic curves $(h,u): \Si \to \RS\times M$, $\Si=\CP-\{(z^{\pm}_k)\}$ 
with two or more punctures, where it remains to understand the map $u$. Note that by proposition 2.2 the perturbed Cauchy-Riemann 
equation for $u$ depends on the $S^1$-component $h_2=h_2^0+t_0$ of the map $h$. Starting with the case of two punctures, we make precise the 
well-known connection between symplectic Floer homology and symplectic field theory for Hamiltonian mapping tori. \\
\\
\begin{proposition} The $\Ju^H$-holomorphic cylinders connecting the $R^H$-orbits $(x^+,T)$ and $(x^-,T)$ in $\RS\times M$ correspond to 
the Floer connecting orbits in $M$ between the one-periodic orbits $x^+(T\cdot)$ and $x^-(T\cdot)$ of the Hamiltonian $H_T(t,\cdot)=T\cdot H(Tt,\cdot)$ and the family $J_T(t,\cdot)=J(Tt,\cdot)$ of $\omega$-compatible almost complex structures. \end{proposition} 
\noindent
{\it Proof:} When $n=2$, i.e., $\uz=(z^-,z^+)$, we find an automorphism $\varphi\in\Aut(\CP)$ with $\varphi(z^-)=0$, 
$\varphi(z^+)=\infty$. Since in the moduli space two elements are considered equal when they agree up to an automorphism of the domain, we can assume 
that $\uz=(0,\infty)$. It follows from lemma 2.3. that $h: \CP-\{0,\infty\} \cong \RS \to \RS$ is of the form 
\begin{equation*} h(s,t)=(Ts+s_0,Tt+t_0) \end{equation*} 
with $T=T^++T^-$. 
We can assume that $h$ is given by $h(s,t)=(Ts,Tt)$ after composing with the automorphism $\varphi(s,t)=(s-s_0/T,t-t_0/T)$ 
of $\RS$. Now the claim follows from the fact that the Cauchy-Riemann equation for $u:\RS\to M$ reads as 
\begin{equation*} \CR_{J,H}u \cdot \del_s = \del_s u + J(Tt,u) \cdot (\del_t u + T\cdot X^H(Tt,u)) = 0, \end{equation*} 
with $T\cdot X^H = X^{T\cdot H}$. $\qed$ \\

\subsection{How to achieve transversality with $S^1$-symmetry}
For understanding the curves with more than two punctures, observe that in these cases the underlying punctured Riemann spheres $\Si$ 
are stable, so that every automorphism $\varphi$ of $\Si$ preserving the ordering of the punctures is the identity. While this implies that different maps $h=h^0+(s_0,t_0)$ give different elements in the moduli space, the main problem is that the solutions for $u$ moreover depend on the $S^1$-component $h_2=h_2^0+t_0$ of the chosen map $h$, that is, the $S^1$-parameter $t_0$. Instead of studying how the solution spaces for $u$ vary with $t_0\in S^1$, it is natural to restrict to special situations when the solution spaces are $t_0$-independent. Moreover, when this can be arranged in a way that all asymptotic orbits are nondegenerate and we can achieve transversality for the moduli spaces, we can use the resulting $S^1$-symmetry on the moduli spaces to show that they do not contribute to the algebraic invariants in rational symplectic field theory. It is easily seen that the Cauchy-Riemann equation is independent of $t_0\in S^1$ when both the family of almost complex structures $J(t,\cdot)$ and the Hamiltonian $H(t,\cdot)$ are independent of $t\in S^1$. Hence for the following we will always assume that 
\begin{equation*} J(t,\cdot) \equiv J,\;\;\; H(t,\cdot) \equiv H. \end{equation*}
\\
{\it Nondegeneracy of the periodic orbits:} \\

It is well-known from symplectic Floer homology that we can achieve that all one-periodic orbits $(x,1)\in P(S^1\times M, H)$ are nondegenerate by choosing $H$ to be a time-independent Morse function $H: M\to\IR$ with a sufficiently small $C^2$-norm, so that, in particular, the only one-periodic orbits of $H$ are the critical points of $H$. While this sounds promising to solve the first of our two problems, note that in contrast to symplectic Floer homology we do not only study curves which are asymptotically cylindrical to one-periodic orbits $(x,1)$ but allow periodic orbits $(x,T)$ of arbitrary period $T\in\IN$. Now the problem is that the $T$-periodic orbits of $H$ are in natural correspondence with one-periodic orbits of the Hamiltonian $T\cdot H$, while $T\cdot H$ need no longer be $C^2$-small enough. In order to solve this problem, we fix a maximal period $T=2^N$ and replace the original Hamiltonian $H$ by $H/2^N$, so that all orbits up to the maximal period $2^N$ are nondegenerate, in particular, critical points of $H/2^N$, i.e., of $H$. \\
\\
{\it Transversality for the Cauchy-Riemann operator:} \\

So it remains the problem of transversality. Although the definition of the algebraic invariants of symplectic field theory suggests that all we have to do is  
counting true $\Ju^H$-holomorphic curves in $\RS\times M$, it is implicit in the definition of all pseudoholomorphic curve theories that before counting 
the geometric data has to be perturbed in such a way that the Cauchy-Riemann operator becomes transversal to the zero section in a suitable Banach space bundle 
over a suitable Banach manifold of maps. It is the main problem of symplectic field theory, as well as Gromov-Witten theory and symplectic Floer 
homology for general symplectic manifolds, that transversality for all moduli spaces cannot be achieved even for generic choices for $\Ju^H$. In fact the problem already occurs for the trivial curves, i.e., trivial examples of curves in symplectic field theory, see [F]. In order to solve these problems virtual moduli cycle techniques were invented, see [LiuT], [LT], [FO]; furthermore they were the starting point for the polyfold project by Hofer, Wysocki and Zehnder, see [H] and the references therein. In order to solve the transversality problem in our $S^1$-symmetric special case, we combine the approach in [CM1] for achieving transversality in Gromov-Witten theory with the well-known connection between symplectic Floer homology and Morse homology in [SZ] as follows. \\
\\
{\it Case of domain-stable curves ($n\geq 3$).} It is well-known, see e.g. [Sch], that transversality in Floer homology and Gromov-Witten theory can be achieved by allowing the almost complex 
structure on the symplectic manifold $(M,\omega)$ to depend on points on the punctured Riemann surface underlying the holomorphic curves, i.e., 
introducing domain-dependent almost complex structures. In this paper we fix the $S^1$-independent almost complex structure $J$ and introduce 
domain-dependent Hamiltonian perturbations $H$, which however are still $S^1$-independent, where we show in section 4 that the resulting class of domain-dependent cylindrical almost complex structures $\Ju^H$ on $\RS\times M$ is still large enough to achieve transversality for all moduli spaces of curves with three or more punctures. Here we let $H$ rather than $J$ depend on the underlying punctured spheres, so that we achieve transversality also for the trivial curves, i.e., the branched covers of trivial cylinders. Note that in order to make the latter transversal, it is clearly necessary to make the stable Hamiltonian structure on $S^1\times M$ domain-dependent. In order to make the choices for the domain-dependent Hamiltonian perturbations $H$ compatible with gluing of curves in symplectic field theory, the perturbations must vary smoothly with the position of the punctures $\uz=(z^{\pm}_1,...,z^{\pm}_{n^{\pm}})$, 
\begin{equation*} H=H_{\uz}: (\CP-\{z^{\pm}_1,...,z^{\pm}_{n^{\pm}}\})\times M\to\IR. \end{equation*} 
In order to guarantee that finite energy solutions are still asymptotically cylindrical over periodic orbits of 
the original domain-independent Hamiltonian $H$, we require that $H_{\uz}$ agrees with $H$ over the cylindrical neighborhoods of the punctures. 
Furthermore, in order to ensure that the automorphism group of $\CP$ still acts on the moduli space, they must satisfy $H_{\varphi(\uz)}= \varphi_* H_{\uz} = H_{\uz}\circ \varphi^{-1}$. When the number of punctures is greater or equal than three, i.e., the punctured Riemann sphere is stable, it follows that $H_{\uz}$ should depend only on the class $[\uz]\in\IM_{0,n}$ in the moduli space of $n$-punctured Riemann spheres. \\
\\
{\it Outline of the construction of domain-dependent Hamiltonians.} For the construction of such domain-dependent structures we follow the ideas in [CM1], where for precise definitions we refer to the upcoming section on domain-dependent Hamiltonian perturbations. For $n\geq 3$ denote by $\IM_{0,n}$ the moduli space of stable genus zero curves modelled over the $n$-labelled tree with one vertex, i.e. the moduli space of Riemann spheres with $n$ marked points. Taking the union of all moduli spaces of stable nodal curves modelled over $n$-labelled trees, we obtain the Deligne-Mumford space $\DM_{0,n} = \coprod_T \IM_T$ which, equipped with the Gromov topology, provides the compactification of 
the moduli space $\IM_{0,n}$. It is a crucial observation that we have a canonical projection $\pi: \DM_{0,n+1} \to \DM_{0,n}$ by forgetting the $(n+1)^{\st}$ marked point and stabilizing. Note that to any nodal curve $\uz$ we can naturally associate a nodal Riemann surface 
$\Sigma_{\uz} = \coprod_{\alpha\in T} S_{\alpha}/\{z_{\alpha\beta}\sim z_{\beta\alpha}\}$ with punctures $(z_k)$,  
obtained by gluing a collection of Riemann spheres $S_{\alpha} \cong \CP$ at the connecting nodes $z_{\alpha\beta} \in \CP$. It then follows that the map $\pi: \DM_{0,n+1} \to \DM_{0,n}$ is holomorphic and the fibre $\pi^{-1}([\uz])$ is naturally biholomorphic to $\Sigma_{\uz}$. We then choose for every $n\geq 3$ smooth maps $H^{(n)}: \DM_{0,n+1} \to C^{\infty}(M)$ and for $[\uz] \in \DM_{0,n}$ 
then define $H_{\uz}$ to be the restriction of $H^{(n)}$ to the fibre $\pi^{-1}([\uz]) \cong \Sigma_{\uz}$. 
In particular, for $\uz \in \IM_{0,n}\subset \DM_{0,n}$ we get from $\Sigma_{\uz} \cong \CP$ a map 
\begin{equation*}
    H_{\uz} = H^{(n)}|_{\pi^{-1}([\uz])}: \CP \to C^{\infty}(M)\,.
\end{equation*}
Assuming that $H^{(2)},...,H^{(n-1)}$ are already chosen, the compatibility for the domain-dependent Hamiltonians under gluing of the underlying Riemann surfaces is ensured by specifying $H^{(n)}$ on the boundary $\del\IM_{0,n+1}=\DM_{0,n+1}-\IM_{0,n+1}$ using $H^{(2)},...,H^{(n-1)}$. For this observe that $\del\IM_{0,n+1}$ consists of the fibres $\pi^{-1}([\uz]) \cong \Sigma_{\uz}$ over $[\uz]\in\del\IM_{0,n} = \DM_{0,n} - \IM_{0,n}$ and of the punctures $z_1,...,z_n \in \CP = \Sigma_{\uz}$ in the fibres over $[\uz]\in\IM_{0,n}$. In order to see that we can indeed define $H^{(n)}$ inductively we crucially use that there are no holomorphic planes and spheres. Assuming we have determined $H^{(n)}$ for $n\geq 2$, we organize all maps into a map 
\begin{equation*} H: \coprod_n \IM_{0,n+1} \to C^{\infty}(M). \end{equation*} 
Note that for $n=2$ the space $\IM_{0,n+1}$ just consists of 
a single point. \\
\\
{\it Case of cylinders ($n=2$).} For curves with two or less punctures, the compatibility with the action of $\Aut(\CP)$ implies that $H_{\uz}$ must be {\it in}dependent 
of points on the domain, i.e., just a function on $M$. On the other hand it is known from symplectic Floer homology that for fixed almost complex structure $J$ it is important to let the Hamiltonian explicitly be $S^1$-dependent to have transversality for generic choices, which seems to destroy our hopes for computing the symplectic field theory of $\RS\times M$ with $S^1$-independent $H$ and $J$. To overcome this problem, we remind ourselves that we already assume $H$ to be so small such that all one-period orbits are nondegenerate, in particular, critical points of $H$. Furthermore by proposition 2.4 we know that the $\Ju^H$-holomorphic cylinders naturally correspond to Floer connecting orbits. The important observation is now to by choosing $H$ with small $C^2$-norm, e.g. by rescaling, we can achieve that all Floer trajectories $u$ are indeed Morse trajectories, i.e., gradient flow lines $u(s,t)\equiv u(s)$ of $H$ between the critical points $x^-$ and $x^+$ with respect to the metric $\omega(\cdot,J\cdot)$ on $M$. When the pair $(H,\omega(\cdot,J\cdot))$ is Morse-Smale, the linearization $F_u$ of the gradient flow operator is surjective, and it is shown in [SZ] that this indeed suffices to show that the linearization $D_u$ of the Cauchy-Riemann operator is surjective as well. More precisely, we use the following lemma, which is proven in [SZ]. \\
\\
\begin{lemma} Let $(H,J)$ be a pair of a Hamiltonian $H$ and an almost complex structure $J$ on a closed symplectic manifold with 
$\<[\omega],\pi_2(M)\>=0$ so that $(H,\omega(\cdot,J\cdot))$ is Morse-Smale. Then the following holds:
\begin{itemize}
\item {\it If $\tau>0$ is sufficiently small, all finite energy solutions $u:\RS \to M$ of $\CR_{J,\tau H}u 
      =\partial_s u + J(u)(\partial_t u + X^{\tau H}(u)) = 0$ are independent of $t\in S^1$.} 
\item {\it In this case, the linearization $D^{\tau}_u$ of $\CR_{J,\tau H}$ is onto 
      at any solution $u:\RS\to M$.} \\
\end{itemize} \end{lemma}

Recall that we fixed a maximal period $T=2^N$ and let $P(H/2^N,\leq 2^N)$ denote the set of periodic orbits of the Reeb vector field $R^{H/2^N}$ 
for the Hamiltonian $H/2^N$ with period less or equal than $2^N$. We collect our results about moduli spaces of holomorphic curves in 
$\RS\times M$ in the following \\
\\
\begin{theorem} Let $(M,\omega)$ be a closed symplectic manifold, which is symplectically aspherical, equipped with 
a $\omega$-compatible almost complex structure $J$ and $H: M\to\IR$ so that lemma 2.5 is satisfied with $\tau=1$. Further assume that for any ordered 
set of punctures $\uz=(z^{\pm}_1,...,z^{\pm}_{n^{\pm}})$ containing three or more points we have constructed a domain-dependent Hamiltonian perturbation 
$H_{\uz}: (\CP-\{\uz\})\times M\to \IR$ of $H$ with the properties outlined above. Then, depending on the 
number of punctures $n$ we have the following result about the moduli spaces of $\Ju^H$-holomorphic curves in $\RS\times M$: 
$ $\\
\begin{itemize}
\item $n=0$: All holomorphic spheres are constant.
\item $n=1$: Holomorphic planes do not exist.
\item $n=2$: For $T\leq 2^N$ the automorphism group $\Aut(\CP)$ acts on the parametrized moduli space 
      $\IM^0(S^1\times M,(x^+,T),(x^-,T),\Ju^{H/2^N})$ of holomorphic cylinders with constant finite isotropy group $\IZ/T\IZ$ and the quotient
      can be naturally identified with the space of gradient flow lines of $H$ with respect to the metric $\omega(\cdot,J\cdot)$ on $M$ 
      between the critical points $x^+$ and $x^-$. 
\item $n \geq 3$: For $P^+,P^- \subset P(H/2^N,\leq 2^N)$ the action of $\Aut(\CP)$ on the parametrized moduli 
      space is free. There still remains a free $S^1$-action on the moduli space after dividing out the $\IR$-translation, where the quotient
      is given by 
      \begin{equation*} \{(u,\uz): u: \CP-\{\uz\} \to M: (*1),(*2)\}/\Aut(\CP)
      \end{equation*}
      {\it with}
      \begin{eqnarray*}
      &(*1):& du + X^{H/2^N}_{\uz}(z,u)\otimes dh_2^0 + J(u) \cdot (du + X^{H/2^N}_{\uz}(z,u)\otimes dh_2^0)\cdot i = 0\,,\\    
      &(*2):& u\circ\psi^{\pm}_k(s,t)\stackrel{s\to\pm\infty}{\longrightarrow} x^{\pm}_k.
      \end{eqnarray*}
\end{itemize} \end{theorem}
\noindent
{\it Proof:} Observe that all statements rely on proposition 2.2 and lemma 2.3. For $n=2$ we additionally use proposition 2.4 and lemma 2.5 
and remark that the critical points and gradient flow lines of $H/2^N$ are naturally identified with those of $H$. 
For the statement about the isotropy groups observe that for $h(s,t)=(Ts,Tt)$ and $u(s,t)=u(s)$ we have  
\begin{equation*} (h,u)=(h\circ\varphi,u\circ\varphi) \Leftrightarrow \varphi(s,t) = (s,t+\frac{k}{T}),\;k\in\IZ/T\IZ. \end{equation*}
For the case $n\geq 3$ observe that the action of $\Aut(\CP)$ is already free on the underlying set of punctures. $\qed$ \\  
 
\section{Domain-dependent Hamiltonians} 
 
Based on the ideas in [CM1] for achieving transversality in Gromov-Witten theory, 
we describe in this section a method to define domain-dependent Hamiltonian perturbations. 
In the following we drop the superscript for the punctures, $\uz=(z_k)$, since for the 
assignment of Hamiltonians we do not distinguish between positive and negative punctures. \\

\subsection{Deligne-Mumford space}

We start with the following definition. \\

\begin{definition} A $n$-labelled tree is a triple $(T,E,\Lambda)$, where $(T,E)$ is a tree with the set of 
vertices $T$ and the 
edge relation $E \subset T\times T$. The set $\Lambda = (\Lambda_{\alpha})_{\alpha\in T}$ is a partition of the index 
set $I=\{1,...,n\}=\bigcup \Lambda_{\alpha}$. We write $\alpha E\beta$ if $(\alpha,\beta)\in E$. \end{definition}

A tree is called {\it stable} if for each $\alpha\in T$ we have $n_{\alpha} = \sharp \Lambda_{\alpha} + 
\sharp\{\beta:\alpha E \beta\} \geq 3$. For $n\geq 3$ a $n$-labelled tree can be stabilized in a canonical way, see [CM1], [MDSa], where one  first deletes vertices $\alpha$ with $n_{\alpha}<3$ to obtain $\st(T) \subset T$ and then modifies $E,\Lambda$ in the obvious way. \\

\begin{definition} A nodal curve of genus zero modelled over $T=(T,E,\Lambda)$ is a tuple $\uz = 
((z_{\alpha\beta})_{\alpha E\beta},(z_k))$ of special points $z_{\alpha\beta}, z_k \in \CP$ such that for each $\alpha\in T$ the special points 
in $Z_{\alpha} = \{z_{\alpha\beta}:\alpha E\beta\}\cup\{z_k: k\in\Lambda_{\alpha}\}$ are pairwise distinct. \end{definition}

To any nodal curve $\uz$ we can naturally associate a nodal Riemann surface 
$\Sigma_{\uz} = \coprod_{\alpha\in T} S_{\alpha}/\{z_{\alpha\beta}\sim z_{\beta\alpha}\}$ with punctures $(z_k)$,  
obtained by gluing a collection of Riemann spheres $S_{\alpha} \cong \CP$ at the points $z_{\alpha\beta} \in \CP$. 
A nodal curve $\uz$ is called {\it stable} if the underlying tree is stable, i.e., every sphere $S_{\alpha}$ carries at 
least three special points. Stabilization of trees immediately leads to a canonical stabilization $\uz \to \st(\uz)$
of the corresponding nodal curve. \\

Denote by $\widetilde{\IM}_T \subset (\CP)^E \times (\CP)^n$ the space of all nodal curves (of genus zero) modelled 
over the tree $T=(T,E,\Lambda)$. An isomorphism between nodal curves $\uz,\uz'$ modelled over the same 
tree is a tuple $\phi = (\phi_{\alpha})_{\alpha\in T}$ with $\phi_{\alpha} \in \Aut(\CP)$ so that 
$\phi(\uz) = \uz'$, i.e., $z'_{\alpha\beta} = \phi_{\alpha}(z_{\alpha\beta})$ and $z'_k 
= \phi_{\alpha}(z_k)$ if $k\in\Lambda_{\alpha}$. Observe that $\phi$ induces a biholomorphism $\phi: 
\Sigma_{\uz} \to \Sigma_{\uz'}$. Let $G_T$ denote the group of biholomorphisms. For stable $T$ the action 
of $G_T$ on $\widetilde{\IM}_T$ is free and the quotient $\IM_T=\widetilde{\IM}_T/G_T$ is a (finite-dimensional) 
complex manifold. \\

\begin{definition} For $n\geq 3$ denote by $\IM_{0,n}$ the moduli space of stable genus zero curves 
modelled over the $n$-labelled tree with 
one vertex, i.e, the moduli space of Riemann spheres with $n$ marked points. Taking the union of all moduli 
spaces of stable nodal curves modelled over $n$-labelled trees, we obtain the Deligne-Mumford space 
\begin{equation*} \DM_{0,n} = \coprod_T \IM_T, \end{equation*} 
which, equipped with the Gromov topology, provides the compactification of 
the moduli space $\IM_{0,n}$ of punctured Riemann spheres. \end{definition}

By a result of Knudsen (see [CM1], theorem 2.1) the Deligne-Mumford space $\DM_{0,n}$ carries the 
structure of a compact complex manifold of complex dimension $n-3$. For each stable $n$-labelled tree $T$ the 
space $\IM_T \subset \DM_{0,n}$ is a complex submanifold, where any $\IM_T \neq \IM_{0,n}$ is 
of complex codimension at least one in $\DM_{0,n}$. \\

It is a crucial observation that we have a canonical projection $\pi: \DM_{0,n+1} \to \DM_{0,n}$ by forgetting the $(n+1)^{\st}$ 
marked point and stabilizing. The map $\pi$ is holomorphic and the fibre $\pi^{-1}([\uz])$ is naturally 
biholomorphic to $\Sigma_{\uz}$. Moreover, for $[\uz] \in \DM_{0,n}$, every component $S_{\alpha} \subset \Sigma_{\uz}$ 
is an embedded holomorphic sphere in $\DM_{0,n+1}$. Note that $\IM_{0,n+1} \varsubsetneqq \pi^{-1}(\IM_{0,n})$ 
as $\pi^{-1}([\uz]) \cap \IM_{0,n+1} = \CP - \{(z_k)\}$ for $[\uz]\in\IM_{0,n}$. \\

\subsection{Definition of coherent Hamiltonian perturbations}

With this we are now ready to describe the algorithm how to 
find domain-dependent Hamiltonians $H_{\uz}$ on $M$. \\ 

For $n=2$ let $H^{(2)}:M\to\IR$ be the domain-{\it in}dependent Hamiltonian from theorem 2.6, i.e., such that with the fixed 
almost complex structure $J$ on $M$ lemma 2.5 is satisfied with $\tau=1$. \\
  
For $n\geq 3$ we choose smooth maps $H^{(n)}: \DM_{0,n+1} \to C^{\infty}(M)$. For $[\uz] \in \DM_{0,n}$ 
we then define $H_{\uz}$ to be the restriction of $H^{(n)}$ to the fibre $\pi^{-1}([\uz]) \cong \Sigma_{\uz}$. 
In particular, for $\uz \in \IM_{0,n}\subset \DM_{0,n}$ we get from $\Sigma_{\uz} \cong \CP$ a map 
\begin{equation*}
    H_{\uz} = H^{(n)}|_{\pi^{-1}([\uz])}: \CP \to C^{\infty}(M)\,,
\end{equation*}
where the biholomorphism $\Sigma_{\uz}\cong\CP$ is fixed by requiring that $(z_1,z_2,z_3)$ are mapped to $(0,1,\infty)$. 
Further let $d_{\uz}=\inf\{d(z_k,z_l): 1\leq k<l \leq n\}$ denote the minimal distance between two marked points with 
respect to the Fubini-Study metric on $\CP$, let $D_{\uz}(z)$ be the ball of radius $d_{\uz}/2$ around 
$z\in\CP$ and set $N_{\uz} = D_{\uz}(z_1)\cup...\cup D_{\uz}(z_n)$. Then we choose $H^{(n)}$ so that 
$H_{\uz}$ agrees with $H^{(2)}$ on $N_{\uz}$. \\

The gluing compatibility is ensured by specifying $H^{(n)}$ on the boundary $\del\IM_{0,n+1}=\DM_{0,n+1}-\IM_{0,n+1}$, which 
consists of the fibres $\pi^{-1}([\uz]) = \Sigma_{\uz}$ over $[\uz]\in\del\IM_{0,n} = \DM_{0,n} - \IM_{0,n}$ and the 
points $z_1,...,z_n \in \CP = \Sigma_{\uz}$ in the fibres over $[\uz]\in\IM_{0,n}$. \\

Note that we have already set 
$H_{\uz}(z_k) = H^{(2)}$. For $[\uz] \in \del\IM_{0,n}=\DM_{0,n} - \IM_{0,n}$ we have 
$H_{\uz} = H^{(n)}|_{\pi^{-1}([\uz])}: \Sigma_{\uz} \to C^{\infty}(M)$ with   
$\Sigma_{\uz} = \coprod S_{\alpha}/\sim$ and $\sharp T\geq 2$. As before let $Z_{\alpha}=\{z^{\alpha}_1,...,
z^{\alpha}_{n_{\alpha}}\}$ denote the set of special points on $S_{\alpha}$. Then we want that 
\begin{equation*} H_{\uz}|_{S_{\alpha}} = H_{\uz^{\alpha}} \end{equation*} 
for $\uz^{\alpha} = (z^{\alpha}_k)$. \\

Since $n_{\alpha}=\sharp Z_{\alpha}<n$, this requirement implies 
that a choice for the map $H^{(n)}: \DM_{0,n+1}\to C^{\infty}(M)$ also fixes the maps $H^{(n')}: \DM_{0,n'+1} 
\to C^{\infty}(M)$ for $n'<n$. \\ 

If $H^{(k)}: \DM_{0,k+1}\to C^{\infty}(M)$, $k=2,...,n-1$ are compatible in the above sense we call them 
coherent. We show how to find $H^{(n)}: \DM_{0,n+1}\to C^{\infty}(M)$ so that $H^{(2)},...,H^{(n)}$ are 
coherent. \\
\\
Let $[\uz]\in\del\IM_{0,n}$ with $\Sigma_{\uz}=\coprod S_{\alpha}/\sim$. Under the assumption 
that $H_{\uz^{\alpha}}$ was chosen to agree with $H^{(2)}$ on the neighborhood $N_{\uz^{\alpha}}$ of the special 
points it follows that all $H_{\uz^{\alpha}}$ fit together to a smooth assignment $H_{\uz}: \Sigma_{\uz} 
\to C^{\infty}(M)$. Let $T = (T,E,\Lambda)$ be the tree underlying $\uz$. Then it follows by the same arguments 
that the maps $H^{(n_{\alpha})}$ fit together to a smooth map $H^T: \pi^{-1}(\DM_T) \to C^{\infty}(M)$. Now let 
$\tau: T \to T'$ be a surjective tree homomorphism with $\sharp T' \geq 2$. Then $\DM_T \subset \DM_{T'}$ and it 
follows from the compatibility of $H^{(2)},...,H^{(n-1)}$ that $H^T$ and $H^{T'}$ agree on $\pi^{-1}(\DM_T)$. Hence 
we get a unique assigment on $\del \IM_{0,n+1} = \pi^{-1}(\coprod \{\IM_T: \sharp T \geq 2\})$. \\
\\
After having specified the map $H^{(n)}:\DM_{0,n+1}\to C^{\infty}(M)$ on the boundary $\del\IM_{0,n+1}$, we 
choose $H^{(n)}$ in the interior $\IM_{0,n+1}$ so that $H^{(n)}$ is smooth (on the compactification $\DM_{0,n+1}$) and 
$H^{(n)}$ agrees with $H^{(2)}$ on $N_{\uz} \subset \pi^{-1}([\uz])$ for all $[\uz]\in\IM_{0,n}$.  \\

Assuming we have determined $H^{(n)}$ for $n\geq 2$, we organize all maps into a map 
\begin{equation*} H: \coprod_n \IM_{0,n+1} \to C^{\infty}(M). \end{equation*} 

Note that for $n=2$ the space $\IM_{0,n+1}$ just consists of 
a single point. A map $H$ as above, i.e., for which all restrictions 
$H^{(n)}: \IM_{0,n+1}\to C^{\infty}(M)$, $n\in\IN$ are coherent, is again called coherent. \\

Together with the almost complex structure $J$ recall that this defines a domain-dependent cylindrical almost complex 
structure $\Ju^H$ on $\RS\times M$, 
\begin{equation*} \Ju^H: \coprod_n \IM_{0,n+1} \to \J_{\cyl}(\RS\times M). \end{equation*} 
With this generalized notion of cylindrical almost complex structure we call, according to theorem 2.6, a map 
$F=(h,u): \CP-\{\uz\}\to \RS\times M$ $\Ju^H$-holomorphic when it satisfies the domain-dependent Cauchy-Riemann equation 
\begin{equation*} \CR_{\Ju} (h,u) = d(h,u)+ \Ju^{H}_{\uz}(z,h,u)\cdot d(h,u) \cdot i= 0, \end{equation*}
which by proposition 2.2 is equivalent to the set of equations $\CR h = 0$ and 
\begin{equation*}
 \CR_{J,H}(u)= du + X^H_{\uz}(z,u)\otimes dh_2^0 + J(u) \cdot (du + X^H_{\uz}(z,u)\otimes dh_2^0)\cdot i = 0 
\end{equation*}
with $X^H_{\uz}(z,\cdot)$ denoting the symplectic gradient of $H_{\uz}(z,\cdot): M\to\IR$. \\
 
Since $H_{\uz}(z,\cdot)$ agrees with the Hamiltonian $H^{(2)}:M\to\IR$ near the punctures, it follows that 
any finite-energy solution of the modified perturbed Cauchy-Riemann equation again converges to a periodic 
orbit of the Hamiltonian flow of $H^{(2)}$ as long as all possible asymptotic orbits are nondegenerate.    
Observe that it follows from the definition of $H_{\uz}$ that the group of Moebius transformations still acts on the resulting moduli 
space of parametrized curves. We show in the section on transversality that for any given almost complex structure $J$ on $M$ we can find Hamiltonian perturbations $H: \coprod_n \IM_{0,n+1} \to C^{\infty}(M)$, so that all moduli spaces $\IM^0(S^1\times M;P^+,P^-;\Ju^{H/2^N})$ are 
cut out transversally simultaneously for all maximal periods $2^N$, $N\in\IN$. \\

\subsection{Compatibility with SFT compactness}

It remains to show that the notion of coherent cylindrical almost complex structures $\Ju^H$ is actually 
compatible with Gromov convergence of $\Ju^H$-holomorphic curves in $\RS\times M$. \\
\\
\begin{definition} A $\Ju^H$-holomorphic level $\ell$ map $(h,u,\uz)$ consists of the following data:
\begin{itemize}
\item A nodal curve $\uz =\coprod S_{\alpha}/\sim \in \DM_{0,n}$ and a labeling $\sigma: T \to \{1,...,\ell\}$,
      called levels, such that two components $\alpha,\beta \in T$ with $\alpha E\beta$ have levels 
      differing by at most one.
\item $\Ju^H$-holomorphic maps $F_{\alpha}=(h_{\alpha},u_{\alpha}): S_{\alpha} \to \RS\times M$ (satisfying  
      $d(h_{\alpha},u_{\alpha}) + \Ju^H_{\uz^{\alpha}}(z,h_{\alpha},u_{\alpha}) \cdot d(h_{\alpha},u_{\alpha})\cdot i = 0$) 
      with the following behaviour at the nodes. \\
      If $\sigma(\alpha) = \sigma(\beta) +1$ then $z_{\alpha\beta}$ is a negative puncture for 
      $(h_{\alpha},u_{\alpha})$ and $z_{\beta\alpha}$ a positive puncture for $(h_{\beta},u_{\beta})$ and they 
      are asymptotically cylindrical over the same periodic orbit; else, if $\sigma(\alpha) = \sigma(\beta)$, then 
      $(h_{\alpha},u_{\alpha})(z_{\alpha\beta}) = (h_{\beta},u_{\beta})(z_{\beta\alpha})$. 
\end{itemize}  \end{definition} 

With this we can give the definition of Gromov convergence of $\Ju^H$-holomorphic maps. \\
\\
\begin{definition} A sequence of stable $\Ju^H$-holomorphic maps $(h^{\nu},u^{\nu},\uz^{\nu})$ converges to 
a level $\ell$ holomorphic map $(h,u,\uz)$ if for any $\alpha \in T$ ($T$ is the tree underlying $\uz$) 
there exists a sequence of Moebius transformations $\phi^{\nu}_{\alpha} \in \Aut(\CP)$ such that:  
\begin{itemize}
\item for $(h,u)=(h_1,h_2,u)=(h_{1,\alpha},h_{2,\alpha},u_{\alpha})_{\alpha\in T}$ there exist sequences 
$s^{\nu}_i$, $i=1,...,\ell$ with 
\begin{eqnarray*}
  h_1^{\nu}\circ\phi^{\nu}_{\alpha} + s^{\nu}_{\sigma(\alpha)} \stackrel{\nu\to\infty}{\longrightarrow} 
  h_{1,\alpha},\,\,\,
  (h_2^{\nu},u^{\nu})\circ\phi^{\nu}_{\alpha} \stackrel{\nu\to\infty}{\longrightarrow} (h_{2,\alpha},u_{\alpha}) 
\end{eqnarray*}
for all $\alpha\in T$ in $C^{\infty}_{\loc}(\Si)$,
\item for all $k=1,...,n$ we have $(\phi^{\nu}_{\alpha})^{-1}(z^{\nu}_k) \to z_k$ if 
      $k\in\Lambda_{\alpha}$ ($z_k \in S_{\alpha}$),
\item and $(\phi^{\nu}_{\alpha})^{-1} \circ \phi^{\nu}_{\beta} \to z_{\alpha\beta}$ for all 
      $\alpha E \beta$.
\end{itemize} \end{definition}

Note that a level $\ell$ holomorphic map $(h,u,\uz)$ is called stable if for any $l\in\{1,...,\ell\}$ 
there exists $\alpha\in T$ with $\sigma(\alpha)=l$ and $(h_{\alpha},u_{\alpha})$ is not a trivial cylinder and, furthermore, 
if $(h_{\alpha},u_{\alpha})$ is constant then the number of special points $n_{\alpha} = \sharp Z_{\alpha} \geq 3$. 
Although any holomorphic map $(h^{\nu},u^{\nu},\uz^{\nu}) \in \IM^0(S^1\times M;P^+,P^-;\Ju^H)$ with 
$n= \sharp P^++\sharp P^- \geq 3$ is stable, the nodal curve $\uz$ underlying the limit level $\ell$ 
holomorphic map $(h,u,\uz)$ need not be stable. However, we can use the absence of holomorphic planes and 
(non-constant) holomorphic spheres in $\RS\times M$ to prove the following lemma about the boundary of 
$\IM(S^1\times M;P^+,P^-;\Ju^H)/\IR$. \\
\\
\begin{lemma} Assume that the sequence $(h^{\nu},u^{\nu},\uz^{\nu}) \in \IM(S^1\times M;P^+,P^-;\Ju^H)$ Gromov 
converges to the level $\ell$ holomorphic map $(h,u,\uz)$. For the number of special points $n_{\alpha}$ 
on the component $S_{\alpha} \subset \Sigma_{\uz}$ it holds
\begin{itemize}
\item $n_{\alpha} \leq n = \sharp P^++\sharp P^-$ for any $\alpha\in T$,
\item if $n_{\alpha} = n$ for some $\alpha\in T$ then all other components are cylinders, i.e., 
      carry precisely two special points.
\end{itemize} \end{lemma}
\noindent
{\it Proof:} We prove this statement by iteratively letting circles on $\CP$ collapse to 
obtain the nodal surface $\Sigma_{\uz}$. \\
For increasing the maximal number of special points on spherical components on a nodal surface we 
must collapse a special circle with all special points on one hemisphere. Even after  
collapsing further circles to nodes there always remains one component with just one 
special point (a node). Since by $\<[\omega],\pi_2(M)\>=0$ there are no holomorphic planes and 
bubbles this cannot happen, which shows the first part of the statement. 
For the second part observe that collapsing circles with more than one special point on each hemisphere leads to two new 
spherical components which carry strictly less special points than the original one. $\qed$ \\

Let $n\geq 3$. For chosen $H: \coprod_n \IM_{0,n+1} \to C^{\infty}(M)$ recall 
that for stable nodal curves $\uz$ we defined $H_{\uz} = H|_{\pi^{-1}([\uz])}: \Sigma_{\uz} 
\to C^{\infty}(M)$. For general nodal curves $\uz$ (with $n\geq 3$) we can use the stabilization $\uz \to 
\st(\uz)$ and the induced map $\st: \Sigma_{\uz} \to \Sigma_{\st(\uz)}$ to define 
\begin{equation*}
  H_{\uz}(z) := H_{\st(\uz)}(\st(z))\,,\,\,\,z\in\Sigma_{\uz}
\end{equation*}    
(compare [CM1], section 4) with corresponding cylindrical almost complex structure 
$\Ju^H_{\uz}(z) := \Ju^H_{\st(\uz)}(\st(z)) \in \J_{\cyl}(S^1\times M)$. \\
\\
\begin{proposition} A $\Ju^H$-holomorphic level $\ell$ map $(h,u,\uz)$ is 
$\Ju^H_{\uz}$-holomorphic. \end{proposition} 
\noindent
{\it Proof:} If $\uz$ is stable this follows directly from the construction of $\Ju^H$ as the restriction 
of $\Ju^H_{\uz}$ to a component $S_{\alpha} \subset \Sigma_{\uz}$ agrees with $\Ju^H_{\uz^{\alpha}}$ when 
$\uz^{\alpha} = (z^{\alpha}_1,...,z^{\alpha}_{n_{\alpha}})$ denotes the ordered set of special points on $S_{\alpha}$. 
If $\uz$ is not stable the proposition relies on the following two observations. \\
Since there are no spherical components with just one special point all  
special points on stable components of $\Sigma_{\uz}$ are preserved under stabilization, i.e., a node connecting 
a stable component with an unstable one is not removed but becomes a marked point on $\Sigma_{\st(\uz)}$. \\
On the other hand points on a cylindrical component (a tree of cylinders) are mapped 
under stabilization to the node connecting it to a stable component (which then is a marked point for the 
nodal surface $\Sigma_{\st(\uz)}$). Since $\Ju^H_{\st(\uz)}$ near special points agrees with complex structure 
$\Ju^{H^{(2)}}$ chosen for cylinder we have $\Ju^H_{\uz}(z) = \Ju^H_{\st(\uz)}(\st(z)) = \Ju^{H^{(2)}}$ for any 
$z\in\Sigma_{\uz}$ lying on a cylindrical component. $\qed$ \\
\\
In order to show the gluing compatibility we prove the following proposition. \\
\\
\begin{proposition} Let $(h^{\nu},u^{\nu},\uz^{\nu})$ be a sequence of 
$\Ju^H_{\uz^{\nu}}$-holomorphic maps converging to the level $\ell$ map $(h,u,\uz)$. Then $(h,u,\uz)$ is 
$\Ju^H_{\uz}$-holomorphic. \end{proposition}
\noindent
{\it Proof:} Recall from the definition of Gromov convergence that for any $\alpha \in T$ 
(the tree underlying $\uz$) there exists a sequence $\phi^{\nu}_{\alpha} \in \Aut(\CP)$ and 
for any $i\in\{1,...,\ell\}$ sequences $s^{\nu}_i \in \IR$ such that 
$h_1^{\nu}\circ\phi^{\nu}_{\alpha} + s^{\nu}_{\sigma(\alpha)} \to h_{1,\alpha}$ and 
$(h_2^{\nu},u^{\nu})\circ\phi^{\nu}_{\alpha} \to (h_{1,\alpha},u_{\alpha})$. Hence it remains to show that 
\begin{equation*}
  \Ju^H_{\uz^{\nu}} \circ \phi^{\nu}_{\alpha} \to \Ju^H_{\uz} 
\end{equation*}
in $C^{\infty}(S_{\alpha},\J_{\cyl}(S^1\times M))$ as $\nu\to\infty$ for all $\alpha \in T$. \\
  
Since the projection from the compactified moduli space to the Deligne-Mumford space $\DM_{0,n}$ is 
smooth (see theorem 5.6.6 in [MDSa]), it follows from $(h^{\nu},u^{\nu},\uz^{\nu}) \to (h,u,\uz)$ that 
$\uz^{\nu} = \st(\uz^{\nu}) \to \st(\uz)$ in $\DM_{0,n}$. \\
For $\alpha\in\st(T)$ and $z\in S_{\alpha}$ we have $\st(z)=z$ and it follows that 
\begin{equation*}
(\uz^{\nu},\phi^{\nu}_{\alpha}(z)) \to (\st(\uz),z) \in \DM_{0,n+1}\,.
\end{equation*}
Since $\Ju^{H^{(n)}}: \DM_{0,n+1} \to \J_{\cyl}(S^1\times M)$ is continuous, we have 
\begin{equation*}
\Ju^H_{\uz^{\nu}}(\phi^{\nu}_{\alpha}(z)) \to \Ju^H_{\st(\uz)}(z) = \Ju^H_{\uz}(z)
\end{equation*}
in $\J_{\cyl}(S^1\times M)$ for all $z\in S_{\alpha}$. The uniform convergence in all derivatives follows 
by the same argument using the smoothness of $\Ju^{H^{(n)}}$. \\
On the other hand, if $\alpha\notin\st(T)$ and $z\in S_{\alpha}$, then $\st(z)=z_{\beta\alpha} \in \st(\uz)$ if $\alpha E\beta$. 
In $\DM_{0,n+1}$ we have that 
\begin{equation*}
  (\uz^{\nu},\phi^{\nu}_{\alpha}(z)) \to (\uz,z_{\beta\alpha})
\end{equation*}
since $(\phi^{\nu}_{\beta})^{-1}(\phi^{\nu}_{\alpha}(z)) \to z_{\beta\alpha} \in S_{\beta}$ 
and therefore 
\begin{equation*}
  \Ju^H_{\uz^{\nu}}(\phi^{\nu}_{\alpha}(z)) \to \Ju^H_{\st(\uz)}(\st(z)) = \Ju^H_{\uz}(z)\,.\,\,\qed
\end{equation*}     

\section{Transversality}

We follow [BM] for the description of the analytic setup of the underlying Fredholm problem. More 
precisely, we take from [BM] the definition of the Banach space bundle over the Banach manifold of maps, 
which contains the Cauchy-Riemann operator studied above as a smooth section. \\

\subsection{Banach space bundle and Cauchy-Riemann operator}

For a chosen coherent Hamiltonian perturbation $H: \coprod_n \IM_{0,n+1}\to C^{\infty}(M)$ and fixed $N\in\IN$, we 
choose ordered sets of periodic orbits 
\begin{equation*}
P^{\pm} =\{(x^{\pm}_1,T^{\pm}_1),...,(x^{\pm}_{n^{\pm}},T^{\pm}_{n^{\pm}})\} \subset P(H^{(2)}/2^N,\leq 2^N),
\end{equation*} 
where $n=n^++n^-$. Instead of considering $\CP \cong S^2$ with its unique conformal structure, we fix punctures $z_1^{\pm,0},...,z_n^{\pm,0} \in S^2$ and let the complex structure on $\Si =S^2-\{z_1^{\pm,0},...,z_n^{\pm,0}\}$ vary. 
Following the constructions in [BM] we see that the appropriate Banach manifold $\BB^{p,d}(\RS\times M;(x^{\pm}_k,T^{\pm}_k))$ 
for studying the underlying Fredholm problem is given by the product 
\begin{equation*}
 \BB^{p,d}(\RS\times M,(x^{\pm}_k,T^{\pm}_k)) = H^{1,p,d}_{\cst}(\Si,\IC) \times \BB^p(M;(x^{\pm}_k)) \times \IM_{0,n}
\end{equation*}
with $d>0$ and $p>2$, whose factors are defined as follows. \\

The Banach manifold $\BB^p(M;(x^{\pm}_k))$ consists of maps $u \in 
H^{1,p}_{\loc}(\Si,M)$, which converge to the critical points $x^{\pm}_k \in \Crit(H^{(2)})$ as $z\in\Si$ approaches 
the puncture $z_k^{\pm,0}$. More precisely, if we fix  
linear maps $\Theta_k^{\pm}: \IR^{2m} \to T_{x^{\pm}_k}M$, the curves satisfy 
\begin{equation*}
  u \circ \psi^{\pm}_k(s,t) = \exp_{x^{\pm}_k}(\Theta_k^{\pm} \cdot v^{\pm}_k(s,t))
\end{equation*}
for some $v^{\pm}_k \in H^{1,p}(\IR^{\pm}\times S^1,\IR^{2m})$, where $\exp$ denotes the exponential map for the 
metric $\omega(\cdot,J\cdot)$ on $M$. \\ 

The space $H^{1,p,d}_{\cst}(\Si,\IC)$ consists of maps $h \in H^{1,p}_{\loc}(\Si,\IC)$, for which there exist $(s_0^{\pm,k},t_0^{\pm,k}) \in \IR^2\cong \IC$, 
such that for all $k=1,...,n^{\pm}$ the maps 
\begin{eqnarray*}
&&\IR^{\pm}\times S^1 \to \IR^2,\,(s,t) \mapsto ((h \circ \psi^{\pm}_k)(s,t) - (s_0^{\pm,k},t_0^{\pm,k})) \cdot e^{\pm d\cdot s} 
\end{eqnarray*}
are in $H^{1,p}(\IR^{\pm}\times S^1,\IC)$. In other words, $H^{1,p,d}_{\cst}(\Si,\IC)$ consists of 
maps differing asymptotically from a constant one by a function, which converges exponentially fast to zero. \\

Finally $\IM_{0,n}$ denotes, as before, the moduli space of complex structures on the punctured sphere $\Si$, which clearly is 
naturally identified with its originally defined version, the moduli space of Riemann spheres with $n$ punctures. \\
 Here we represent $\IM_{0,n}$ explicitly by finite-dimensional families of (almost) complex structures on $\Si$, 
so that $T_j\IM_{0,n}$ becomes a finite-dimensional subspace of 
\begin{equation*} \{y \in \End(T\Si): yj+jy=0\}. \end{equation*}
Note that in [BM] the authors work with Teichmueller spaces, since the corresponding moduli spaces of complex 
structures, obtained by dividing out the mapping class group, become orbifolds for non-zero genus. \\   

Note that for the identification
\begin{equation*}
 \BB^{p,d}(\RS\times M,(x^{\pm}_k,T^{\pm}_k)) = H^{1,p,d}_{\cst}(\Si,\IC) \times \BB^p(M;(x^{\pm}_k)) \times \IM_{0,n}
\end{equation*}
we identify $\hb\in H^{1,p,d}_{\cst}(\Si,\IC)$ with the map $h: \Si\to\RS$ given by $h=h^0+\hb$, where $h^0$ denotes an arbitrary fixed holomorphic map $h^0: \Si \to \RS\cong \CP-\{0,\infty\}$, so that $z_k^{\pm,0}$ is a pole/zero of order $T_k^{\pm}$. Note that we do not use asymptotic exponential weights (depending on $d \in \IR^+$) for the Banach manifold $\BB^p(M;(x^{\pm}_k))$, since we 
are dealing with nondegenerate asymptotics. \\

Let $H^{1,p}(u^*TM)$ consist of sections $\xi \in H^{1,p}_{\loc}(u^*TM)$, such that
\begin{equation*}
\xi \circ \psi^{\pm}_k(s,t) = (d\exp_{x^{\pm}_k})(\Theta_k^{\pm}\cdot v_k^{\pm}(s,t)) \cdot \Theta_k^{\pm}\xi_k^{\pm,0}(s,t)
\end{equation*}
with $\xi^{\pm,0}_k \in H^{1,p}(\IR^{\pm}\times S^1,\IR^{2m})$ for $k=1,...,n$. Note that here we take the differential 
of $\exp_{x^{\pm}_k}: T_{x^{\pm}_k}M \to M$ at $\Theta_k^{\pm}\cdot v_k^{\pm}(s,t) \in T_{x^{\pm}_k}M$, which maps 
the tangent space to $M$ at $x^{\pm}_k$ to the tangent space to $M$ at 
\begin{equation*} \exp_{x^{\pm}_k}(\Theta_k^{\pm}\cdot v_k^{\pm}(s,t))= u\circ\psi^{\pm}_k(s,t). \end{equation*} 

Then the tangent space to 
$\BB^{p,d}(\RS\times M;(x^{\pm}_k,T^{\pm}_k))$ at $(\hb,u,j)$ is given by
\begin{equation*}
  T_{(\hb,u,j)} \BB^{p,d}(\RS\times M;(x^{\pm}_k,T^{\pm}_k)) =  
  H^{1,p,d}_{\cst}(\Si,\IC) \oplus H^{1,p}(u^*TM) \oplus T_j\IM_{0,n}. 
\end{equation*}
 
Consider the bundle $\Lambda^{0,1}\otimes_{j,J} u^*TM$, 
whose sections are (0,1)-forms on $\Si$ with values in the pullback bundle $u^*TM$ equipped with the complex structure $J$. 
The space $L^p(\Lambda^{0,1}\otimes_{j,J} u^*TM)$ is defined similarly as $H^{1,p}(u^*TM)$: it consists of sections 
$\alpha \in L^p_{\loc}$, which asymptotically satisfy 
\begin{equation*}
(\psi^{\pm}_k)^*\alpha(s,t) \cdot \partial_s = (d\exp_{x^{\pm}_k})(\Theta_k^{\pm}\cdot v_k^{\pm}(s,t)) \cdot 
\Theta_k^{\pm}\alpha^{\pm,0}_k(s,t)
\end{equation*} 
with $\alpha_k^{\pm,0} \in L^p(\IR^{\pm}\times S^1,\IR^{2m})$. \\

Over $\BB^{p,d} = \BB^{p,d}(\RS\times M;(x^{\pm}_k,T^{\pm}_k))$ consider 
the Banach space bundle $\EE^{p,d} \to \BB^{p,d}$ with fibre 
\begin{equation*}
\EE^{p,d}_{\hb,u,j} = L^{p,d}(\Lambda^{0,1}\otimes_{j,i}\IC) \oplus 
L^p(\Lambda^{0,1}\otimes_{j,J} u^*TM).
\end{equation*}

Assume that we have fixed a coherent Hamiltonian perturbation $H: \coprod \IM_{0,n+1} \to C^{\infty}(M)$.  Our convention at the beginning 
of this section, i.e., fixing the punctures on $S^2$ but letting the almost complex structure $j: T\Si\to T\Si$ vary, now 
leads to a dependency $H(j,z)=H^{(n)}(j,z)$ on the complex structure $j$ on $\Si$ and points $z \in \Si$. For the 
following exposition let us assume $N=0$ in order to keep the notation simple. \\
 
The Cauchy-Riemann operator 
\begin{eqnarray*}
 \CR_{\Ju^H}(h,u,j) = \CR_{j,\Ju^H} (h,u) &=& d(h,u) + \Ju^H(j,z,h,u) \cdot d(h,u) \cdot j 
\end{eqnarray*}
is a smooth section in $\EE^{p,d} \to \BB^{p,d}$ and naturally splits, 
\begin{equation*} 
 \CR_{j,\Ju^H}(h,u) = (\CR h,\CR_{J,H}u) \in L^{p,d}(\Lambda^{0,1}\otimes_{j,i}\IC) \oplus 
L^p(\Lambda^{0,1}\otimes_{j,J} u^*TM). 
\end{equation*}

Here $\CR = \CR_{j,i}$ is the standard Cauchy-Riemann operator for maps $h:(\Si,j)\to\RS$ and $\CR_{J,H}$ is the 
perturbed Cauchy-Riemann operator given by 
\begin{equation*}
\CR_{J,H}(u) = du + X^H (j,z,u) \otimes dh_2^0 + J(u)\cdot(du+X^H (j,z,u)\otimes dh_2^0)\cdot j,
\end{equation*}
where again $X^H(j,z,\cdot)$ denotes the symplectic gradient of $H(j,z,\cdot): M\to\IR$. It follows that the linearization $D_{\hb,u,j}$ of $\CR_{\Ju^H}$ at a solution $(\hb,u,j)$ splits, 
\begin{equation*} D_{\hb,u,j}=D_{\hb,u}\oplus D_{j}, \end{equation*} 
with $D_j: T_j \IM_{0,n} \to \EE^{p,d}_{\hb,u,j}$ and 
\begin{eqnarray*}
D_{\hb,u} = \diag(\CR,D_u): &&H^{1,p,d}_{\cst}(\Si,\IC) \oplus H^{1,p}(u^*TM) \\
 &&\to L^{p,d}(\Lambda^{0,1}\otimes_{j,i}\IC) \oplus L^p(\Lambda^{0,1}\otimes_{j,J} u^*TM), 
\end{eqnarray*}
where
\begin{eqnarray*}
 &D_u:& H^{1,p}(u^*TM) \to L^p(\Lambda^{0,1}\otimes_{j,J} u^*TM), \\
 &D_u \xi& = \nabla\xi + J(u) \cdot \nabla\xi \cdot j 
             + \nabla_{\xi}J(u)\cdot du \cdot j \\&& + \nabla_{\xi}X^{H}(j,z,u) \otimes dh_2^0 
             + \nabla_{\xi}\nabla H(j,z,u) \otimes dh_1^0
\end{eqnarray*}
is the linearization of the perturbed Cauchy-Riemann operator $\CR_{J,H}$.

\subsection{Universal moduli space}

Let $\H^{\ell}_n(M;H^{(2)},...,H^{(n-1)})$ denote the Banach manifold consisting of $C^{\ell}$-maps 
$H^{(n)}: \IM_{0,n+1}\to C^{\ell}(M)$, which extend as $C^{\ell}$-maps to $\DM_{0,n+1}$ as induced by $H^{(k)}$, 
$k=2,...,n-1$ and $H^{(n)}(j,\cdot) = H^{(2)}$ on a neighborhood $N_0 \subset \Si$ of the punctures. \\

Note that it is essential to work in the $C^{\ell}$-category since the corresponding space of $C^{\infty}$-structures 
just inherits the structure of a Frechet manifold and we later cannot apply the Sard-Smale theorem. \\
 
The tangent space to $\H^{\ell}=\H^{\ell}_n(M;H^{(2)},...,H^{(n-1)})$ at $H = H^{(n)}$ is given by 
\begin{eqnarray*}
  T_H\H^{\ell}_n(M;H^{(2)},...,H^{(n-1)}) &=& \H^{\ell}_n(M;0,...,0).
\end{eqnarray*}

The universal Cauchy-Riemann operator $\CR_J(\hb,u,j,H) := \CR_{\Ju^H}(h,u,j)$ extends to a smooth section in the 
Banach space bundle $\hat{\EE}^{p,d} \to \BB^{p,d} \times \H^{\ell}$ with fibre 
\begin{equation*}
    \hat{\EE}^{p,d}_{\hb,u,j,H} = \EE^{p,d}_{\hb,u,j} = L^{p,d}(\Lambda^{0,1}\otimes_{j,i}\IC) \oplus 
    L^p(\Lambda^{0,1}\otimes_{j,J} u^*TM).
\end{equation*}      

Letting $\Ju^{H^{(2)}},...,\Ju^{H^{(n-1)}}$ denote the domain-dependent cylindrical almost complex structures on $\RS\times M$ induced by $J$ and $H^{(2)},...,H^{(n-1)}$, we define the universal moduli space $\IM(S^1\times M;P^+,P^-;\Ju^{H^{(2)}},...,\Ju^{H^{(n-1)}})$ as the zero set of the 
universal Cauchy-Riemann operator, 
\begin{eqnarray*} 
 &&\IM(S^1\times M;P^+,P^-;(\Ju^{H^{(k)}})_{k=2}^{n-1}) = \\
 &&\{(\hb,u,j,H)\in 
 \BB^{p,d}\times\H^{\ell}: \CR_J(\hb,u,j,H)=0\}. 
\end{eqnarray*} 

\begin{theorem} For $n\geq 3$ let $H^{(2)},...,H^{(n-1)}$ be fixed. Then for any chosen $(P^+,P^-)$ with $\sharp P^+ + \sharp P^- = n$, 
the universal moduli space $\IM(S^1\times M;P^+,P^-;(\Ju^{H^{(k)}})_{k=2}^{n-1})$ is transversally cut out by the universal Cauchy-Riemann operator 
$\CR_J: \BB^{p,d} \times \H^{\ell} \to \hat{\EE}^{p,d}$ for $d>0$ sufficiently small. In particular, it carries the structure of a $C^{\infty}$-Banach manifold. \end{theorem}

The proof relies on the following two lemmata. \\

\begin{lemma} The operator $\CR: H^{1,p,d}_{\cst}(\Si,\IC) \to L^{p,d}(\Lambda^{0,1}\otimes_{j,i}\IC)$ is onto. \end{lemma}
\noindent
{\it Proof:} Fix a splitting
\begin{equation*}
 H^{1,p,d}_{\cst}(\Si,\IC) = H^{1,p,d}(\Si,\IC) \oplus \Gamma^n
\end{equation*}
where $\Gamma^n\subset C^{\infty}(\Si,\IC)$ is a $2n$-dimensional space of functions storing the constant shifts 
(see [BM]). Given a function $\varphi_d: \Si \to \IR$ 
with $(\varphi_d \circ \psi^{\pm}_k)(s,t) = e^{\pm d\cdot s}$, multiplication with $\varphi_d$ defines 
isomorphisms
\begin{eqnarray*}
H^{1,p,d}(\Si,\IC) &\stackrel{\cong}{\longrightarrow}& H^{1,p}(\Si,\IC), \\
L^{p,d}(T^*\Si\otimes_{i,i}\IC) &\stackrel{\cong}{\longrightarrow}& L^p(T^*\Si\otimes_{i,i}\IC),
\end{eqnarray*}
under which $\CR$ corresponds to a perturbed Cauchy-Riemann operator 
\begin{equation*}
 \CR_d = \CR + S_d: H^{1,p}(\Si,\IC) \to L^p(T^*\Si\otimes_{i,i}\IC).
\end{equation*}
With the asymptotic behaviour of $\varphi_d$ one computes
\begin{equation*}
 S^{\pm,k}_d(t) = (S_d \circ \psi^{\pm}_k)(\pm\infty,t) = \diag(\mp d,\mp d)
\end{equation*}
so that the Conley-Zehnder index for the corresponding paths $\Psi^{\pm,k}: \IR \to \Sp(2m)$ of 
symplectic matrices is $\mp 1$ for $d>0$ sufficiently small. Hence the index of $\CR: 
H^{1,p,d}_{\cst}(\Si,\IC) \to L^{p,d}(T^*\Si\otimes_{i,i}\IC)$ is given by
\begin{equation*}
\ind \CR = \dim \Gamma^n + \ind \CR_d = 2n + \bigl(-n + 1\cdot(2-n)\bigr) = 2,
\end{equation*}
where the sum in the big bracket is the usual index formula for Cauchy-Riemann type operators. On the other hand, it follows from Liouville's theorem that the kernel of $\CR$ consists of the constant functions on $\Si$, so that $\dim \coker \CR = 0$. $\qed$ \\
\\
\begin{lemma} For $n\geq 3$ the linearization $D_{u,H}$ of $\CR_J(u,H) = \CR_{J,H}(u)$ is surjective 
at any $(\hb,u,j,H) \in \IM(S^1\times M;P^+,P^-;(\Ju^{H^{(k)}})_{k=2}^{n-1})$. \end{lemma}
\noindent
{\it Proof:} The operator $D_{u,H}$ is the sum of the linearization $D_u$ of the perturbed 
Cauchy-Riemann operator $\CR_{J,H}$ and the linearization of $\CR_J$ in the $\H^{\ell}$-direction,   
\begin{eqnarray*}
 &D_H:& T_H\H^{\ell} \to L^p(\Lambda^{0,1}\otimes_{j,J} u^*TM),\\ 
 && D_H G = X^G(j,z,u) \otimes dh_2^0 + J(u) X^G(j,z,u)\otimes dh_1^0\,.
\end{eqnarray*}
\\
We show that $D_{u,H}$ is surjective using well-known arguments. Since $D_u$ is Fredholm, the range of $D_{u,H}$ in $L^p(\Lambda^{0,1}\otimes_{j,J} u^*TM)$ is closed, and it suffices to prove that the annihilator of the range of $D_{u,H}$ is trivial. \\

We identify the dual space 
of $L^p(\Lambda^{0,1}\otimes_{j,J} u^*TM)$ with $L^q(\Lambda^{0,1}\otimes_{j,J} u^*TM)$, $1/p+1/q=1$ using the $L^2$-inner 
product on sections in $\Lambda^{0,1}\otimes_{j,J} u^*TM$, which is defined using the standard hyperbolic metric on $(\Si,j)$ 
and the metric $\omega(\cdot, J\cdot)$ on $M$. \\ 

Let $\eta\in\hat{\EE}^{q,d}_{\hb,u,j,H} 
= L^{q,d}(\Lambda^{0,1}\otimes_{j,i}\IC) \oplus L^q(\Lambda^{0,1}\otimes_{j,J} u^*TM)$ such that 
\begin{equation*} \<D_{u,H}\cdot(\xi,G),\eta\>=0 \end{equation*} 
for all $\xi\in H^{1,p}(u^*TM)$ and $G\in T_H\H^{\ell}$. Then surjectivity of $D_{u,H}$ is equivalent to proving $\eta\equiv 0$. \\  

From $\<D_u\xi,\eta\>=0$ for all $\xi\in H^{1,p}(u^*TM)$, we get that $\eta$ is a weak solution of 
the perturbed Cauchy-Riemann equation $D_u^*\eta=0$, where $D_u^*$ is the adjoint of 
$D_u$. By elliptic regularity, it follows that $\eta$ is smooth and hence a strong solution. By unique continuation,  
which is an immediate consequence of the Carleman similarity principle, it follows that $\eta \equiv 0$ whenever 
$\eta$ vanishes identically on an open subset of $\Si$. \\

On the other hand we have 
\begin{eqnarray*} 
 0 \,=\, \<D_H G,\eta\> &=& \int_{\Si} \<J(u) X^G(j,z,u)\otimes dh_1^0 + X^G(j,z,u) \otimes dh_2^0,\eta(z)\>\;dz \\
                        &=& \int_{\Si} \<\nabla G(j,z,u)\otimes dh_1^0 - J(u)\nabla G(j,z,u) \otimes dh_2^0,\eta(z)\>\;dz 
\end{eqnarray*}
for all $G\in T_H\H^{\ell}$. When $z\in\Si$ is not a branch point of the map $h^0: \Si\to\RS$, observe that we can write 
\begin{equation*} \eta(z) = \eta_1(z) \otimes dh^0_1 + \eta_2(z) \otimes dh^0_2 \end{equation*} 
with $\eta_2(z) + J(u)\eta_1(z) = 0$, since $\eta$ is a (0,1)-form. It follows that 
\begin{eqnarray*}
 &&\<\nabla G(j,z,u)\otimes dh_1^0 - J(u)\nabla G(j,z,u) \otimes dh_2^0,\eta(z)\> \\
 &&= \<\nabla G(j,z,u)\otimes dh_1^0 - J(u)\nabla G(j,z,u) \otimes dh_2^0, \eta_1(z) \otimes dh^0_1 - J(u)\eta_1(z) \otimes dh^0_2\> \\ 
 &&= \<\nabla G(j,z,u),\eta_1(z)\>\cdot\|dh_1^0\|^2 \;+\; \<J(u)\nabla G(j,z,u),J(u) \eta_1(z)\> \cdot\|dh_2^0\|^2\\
 &&= \|dh^0\|^2\,\cdot\<\nabla G(j,z,u),\eta_1(z)\> \;=\;\|dh^0\|^2\,\cdot dG(j,z,u)\cdot\eta_1(z),
\end{eqnarray*}
where $\|dh^0\|^2=\|dh_1^0\|^2+\|dh_2^0\|^2$ and $dG(j,z,\cdot)$ denotes the differential of $G(j,z,\cdot): M\to\IR$. \\
 
With this we now prove that $\eta$ vanishes identically on the complement of the set of branch points of $h^0$, which 
by unique continuation implies $\eta=0$. Assume to the contrary that $\eta(z_0)\neq 0$ for some $z_0\in \Si$, which is not a branch point. Since $\alpha$ is a (0,1)-form it follows that 
$\eta_1(z_0)\neq 0$ and we obviously can find $G_0 \in C^{\infty}(M)$ such that 
\begin{equation*} dG_0(u(z_0)) \cdot \eta_1(z_0) > 0. \end{equation*}

Setting $j_0:=j$, let $\varphi\in C^{\infty}(\DM_{0,n+1},[0,1])$ be a smooth cut-off function around 
$(j_0,z_0)\in \IM_{0,n+1}$ 
with $\varphi(j_0,z_0)=1$ and $\varphi(j,z)=0$ for $(j,z) \not\in U(j_0,z_0)$. Here the neighborhood 
$(j_0,z_0)\in U_1(j_0)\times U_2(z_0) = U(j_0,z_0) \subset \DM_{0,n+1}$ is chosen so small that
\begin{equation*} U(j_0,z_0)\cap (\DM_{0,n+1}-\IM_{0,n+1}) =\emptyset,\,\, U_2(z_0)\cap  N_0 = \emptyset, \end{equation*} 
and $dG_0(z,u(z))\cdot\eta_1(z) \geq 0$ for all $z \in U_2(z_0)$. \\

With this define $G: \DM_{0,n+1}\times M \to \IR$ by $G(j,z,p) := 
\varphi(j,z) \cdot G_0(p)$. But this leads to the desired contradiction since we found 
$G \in T_H H^{\ell}= \H^{\ell}_n(M;0,...,0)$ with 
\begin{eqnarray*}
 \<D_H\cdot G,\eta\> = \int_{U_2(z_0)} \|dh^0(z)\|^2\,\cdot\,dG(j,z,u)\cdot \eta_1(z) \;dz > 0. \;\;\qed
\end{eqnarray*} 
$ $\\
{\it Proof of theorem 4.1:} For $n\geq 3$ we must show that the linearization $D_{\hb,u,j,H}$ of the 
universal Cauchy-Riemann operator $\CR_J$ 
is surjective at any \\$(\hb,u,j,H) \in \IM(S^1\times M;P^+,P^-;(\Ju^{H^{(k)}})_{k=2}^{n-1})$. 
Using the splitting $D_{\hb,u,j,H} = D_{\hb,u,H} + D_j$ we show that the first summand 
\begin{eqnarray*}
 D_{\hb,u,H}: &&H^{1,p,d}_{\cst}(\Si,\IC)\oplus T_u\BB^p(M;P^+,P^-)\oplus T_H\H^{\ell} \\
 &&\to L^{p,d}(\Lambda^{0,1}\otimes_{j,i}\IC) \oplus L^p(\Lambda^{0,1}\otimes_{j,J} u^*TM)
\end{eqnarray*}        
is onto. However, since 
\begin{equation*}
D_{\hb,u,H} =  \diag(\CR,D_{u,H}),
\end{equation*}
this follows directly from the surjectivity of $\CR$ and $D_{u,H} = D_u+D_H$. $\qed$\\ 

The importance of the above theorem is that, combined with lemma 2.5, we obtain transversality for all moduli spaces of holomorphic curves in 
$\RS\times M$ asymptotically cylindrical over periodic orbits up to the given maximal period $2^N$. Moreover we can achieve that this 
holds for all maximal periods simultaneously. \\
\\  
\begin{corollary} For $n=2$ and $T\leq 2^N$ the moduli spaces \\$\IM(S^1\times M;(x^+,T),(x^-,T);\Ju^{H/2^N})$ are transversally cut out by 
the Cauchy-Riemann operator for all $N\in\IN$. For $n\geq 3$ we can choose $H^{(n)}\in\H^{\ell}$, simultaneously for all $N\in\IN$, so that the moduli spaces $\IM(S^1\times M;P^+,P^-;\Ju^{H/2^N})$ are transversally cut out by the resulting Cauchy-Riemann operator for all 
$P^+,P^-\subset P(H^{(2)}/2^N,\leq 2^N)$ with $\#P^++\#P^-=n$. \end{corollary}
\noindent
{\it Proof:} For $n=2$ the linear operator 
\begin{equation*}
D_{\hb,u} = \diag(\CR,D_u)
\end{equation*}
is surjective since $D_u$ is onto by lemma 2.5. Indeed, recall that we have chosen the pair $(H^{(2)},J)$ to be regular in the 
sense that $(H^{(2)},\omega(\cdot,J\cdot))$ is Morse-Smale, which 
implies that all pairs $(H^{(2)}/2^N,J)$ for any $N\in\IN$ are again regular, since the stable and unstable manifolds 
are the same. \\

For $n\geq 3$ and $N=0$ the Sard-Smale theorem applied to the map 
\begin{equation*}
\IM(S^1\times M;P^+,P^-;(\Ju^{H^{(k)}})_{k=2}^{n-1}) \to \H^{\ell}_n(M;(H^{(k)})_{k=2}^{n-1}),\,\, 
(\hb,u,j,H)\mapsto H
\end{equation*}
tells us that the set of Hamiltonian perturbations $\H^{\ell}_{\reg}(P^+,P^-) = \\\H^{\ell}_{\reg}(P^+,P^-,0)$, for which 
the moduli space $\IM(S^1\times M;P^+,P^-;\Ju^H)$ is cut out transversally by the Cauchy-Riemann operator 
$\CR_{\Ju^H}$, is of the second Baire category in $\H^{\ell} = \H^{\ell}_n(M;(H^{(k)})_{k=2}^{n-1})$. 
Since there exist just a countable number of tuples $(P^+,P^-)$ with $\sharp P^++\sharp P^-=n$, it follows that 
$\H^{\ell}_{\reg}=\H^{\ell}_{\reg}(0) = \bigcap\{\H^{\ell}_{\reg}(P^+,P^-,0): \sharp P^++\sharp P^-=n\}$ 
is still of the second category. \\

Replacing $H^{(2)},...,H^{(n-1)}$ in the above argumentation by $H^{(2)}/2^N,...,\\H^{(n-1)}/2^N$ for each $N\in\IN$, we
obtain sets of regular structures $\H^{\ell}_{\reg}(N)$, for which the moduli 
spaces $\IM(S^1\times M;P^+,P^-;\Ju^{H/2^N})$ are cut out transversally for all $P^+,P^-\subset P(H^{(2)}/2^N,\leq 2^N)$. 
However, it follows that $\H^{\ell}_{\reg}= \bigcap\{\H^{\ell}_{\reg}(N): N\in\IN\}$ is still of the second category in $\H^{\ell}$. $\qed$  \\

\section{Cobordism}

Since our statements only hold up to a maximal period for the asymptotic orbits, we cannot use the same coherent Hamiltonian 
perturbation to compute the full contact homology. As seen above we must rescale the Hamiltonian for the cylindrical moduli spaces, 
which clearly affects the Hamiltonian perturbations for all punctured spheres. For showing that the graded vector space isomorphism we obtain 
is actually an isomorphism of graded algebras, we construct chain maps between the differential algebras for the different 
coherent Hamiltonian perturbations, which are defined by counting holomorphic curves in an almost complex manifold with 
cylindrical ends.

\subsection{Moduli spaces}

For a given Hamiltonian $H:M\to\IR$ let $\tilde{H}:\IR\times M\to\IR$ be a smooth homotopy with $\tilde{H}(s,\cdot) = H/2$ for $s \leq -1$ and 
$\tilde{H}(s,\cdot)=H$ for $s \geq +1$. Besides that $\tilde{H}$ defines a homotopy of stable Hamiltonian structures 
$(\omega^{\tilde{H}},\lambda^{\tilde{H}})$ with corresponding (constant) symplectic hyperplane bundles $\xi^{\tilde{H}}=TM$ and 
$\IR$-dependent Reeb vector fields $R^{\tilde{H}}(s,t,p) = \del_t + X^{\tilde{H}}(s,t,p)$, it equips $\RS\times M$ with the structure of 
a symplectic manifold with stable cylindrical ends 
\begin{equation*}
((-\infty,-1]\times S^1\times M,\omega^{H/2},\lambda^{H/2})\;\;{\rm and}\;\; ([+1,+\infty)\times S^1\times M,\omega^H,\lambda^H), 
\end{equation*}
where the symplectic structure on the compact, non-cylindrical part 
$(-1,+1)\times S^1\times M$ is given by 
\begin{equation*} \underline{\omega}^{\tilde{H}} = \omega^{\tilde{H}} + ds\wedge dt \end{equation*} 
with $\omega^{\tilde{H}} = \omega + d\tilde{H} \wedge dt$. \\ 
 
Together with the fixed $\omega$-compatible almost complex structure $J$ on $M$, the homotopy $\tilde{H}$ further equips $\RS\times M$ 
with an almost complex structure $\Ju^{\tilde{H}}$ by requiring that it turns $\xi^{\tilde{H}}=TM$ into a complex subbundle with complex 
structure $J$ and 
\begin{equation*} \Ju^{\tilde{H}}\cdot\del_s = R^{\tilde{H}}(s,\cdot) = \del_t + X^{\tilde{H}}(s,\cdot). \end{equation*}
It follows that $(\RS\times M,\Ju^{\tilde{H}})$ is an almost complex manifold with cylindrical ends $((-\infty,-1]\times S^1\times M,\Ju^{H/2})$ 
and $([+1,+\infty)\times S^1\times M,\Ju^H)$. Note that $\Ju^{\tilde{H}}$ is indeed $\underline{\omega}^{\tilde{H}}$-compatible. \\

For our applications we clearly have to replace the Hamiltonian $H:M\to\IR$ by the domain-dependent Hamiltonian perturbation 
$H: \coprod_n \IM_{0,n+1} \times M \to \IR$ from before. It follows that the Hamiltonian homotopy $\tilde{H}$ has to depend explicitly on points 
on the underlying stable punctured spheres, i.e., for the following we consider coherent Hamiltonian homotopies 
\begin{equation*}
 \tilde{H}: \coprod_n \IM_{0,n+1} \times \IR \times M \to \IR,
\end{equation*}
with corresponding domain-dependent almost complex structures 
\begin{equation*} \Ju^{\tilde{H}}: \coprod_n\IM_{0,n+1}\to \J(S^1\times M). \end{equation*}

While it is again clear that the moduli spaces of $\Ju^{\tilde{H}}$-holomorphic curves with more than two punctures come with an $S^1$-symmetry, 
it remains to verify nondegeneracy for the asymptotic orbits and transversality for the curves. Note for the first that we again have to consider 
rescaled versions $\tilde{H}_N: \coprod_n \IM_{0,n+1} \times \IR \times M \to \IR$ with $\tilde{H}_N(s) = \tilde{H}(s/2^N)/2^N$. Since 
$\tilde{H}_N(s) = H/2^{N+1}$ for $s\leq -2^N$ and $\tilde{H}_N(s) = H/2^N$ for $s\geq +2^N$, it is clear that the nondegeneracy holds for all 
asymptotic orbits of period less or equal to $2^N$. \\

While we show below that we can again achieve transversality for all $\Ju^{\tilde{H}}$-holomorphic curves with more than three punctures making 
use of the domain-dependency of the almost complex structure, it remains to guarantee transversality for $\Ju^{\tilde{H}}$-holomorphic 
cylinders. Note that in analogy to proposition 2.4 it follows that all $\Ju^{\tilde{H}}_N$-holomorphic cylinders connecting orbits 
$(x^+,T)$ and $(x^-,T)$ with $T\leq 2^N$ are in natural correspondence to cylinders in $M$ connecting the critical points $x^+,x^-$, which satisfy 
the $\IR$-dependent perturbed Cauchy-Riemann equation 
\begin{equation*}
\CR_{J,H}u \cdot \del_s = \del_s u + J(u) \cdot (\del_t u + T\cdot X^{\tilde{H}}(Ts,u)) = 0.
\end{equation*}
While in general transversality generically only holds for $t$-dependent Hamiltonian homotopies $\tilde{H}$, we can now make use of the following 
natural generalization of lemma 2.5. \\
\\
\begin{lemma} Let $(H,J)$ be a pair of a Hamiltonian $H$ and an almost complex structure $J$ on a closed symplectic manifold with 
$\<[\omega],\pi_2(M)\>=0$ so that $(H,\omega(\cdot,J\cdot))$ is Morse-Smale. Choose $\varphi \in C^{\infty}(\IR,\IR^+)$ with $\varphi(s)=1/2$ 
for $s\leq -1$ and $\varphi(s)=1$ for $s\geq 1$, and let $\tilde{H}: \IR\times M\to \IR$, $\tilde{H}(s,p) = \varphi(s)\cdot H(p)$. 
Then the following holds: 
\begin{itemize}
\item The linearization $\tilde{F}_u$ of $\nabla_{J,\tilde{H}}u = \partial_s u + J(u)X^{\tilde{H}}(s,u)$ is surjective 
      at all solutions.
\item If $\tau>0$ is sufficiently small, all finite energy solutions $u:\RS \to M$ of $\CR_{J,\tilde{H}^{\tau}}u 
      =\partial_s u + J(u)(\partial_t u + X^{H^{\tau}}(s,u)) = 0$ with $\tilde{H}^{\tau}(s,\cdot)=\tau \tilde{H}(\tau s,\cdot)$ 
      are independent of $t\in S^1$. 
\item In this case, the linearization $\tilde{D}_u = \tilde{D}^{\tau}_u$ of $\CR_{J,\tilde{H}^{\tau}}$ is onto 
      at any solution $u:\RS\to M$.
\end{itemize} \end{lemma}
\noindent
{\it Proof:} The proof is a simple generalization of the arguments given in [SZ] and we just show the first statement.  
Let $\tilde{\varphi}: \IR\to\IR^+$ with $\del_s\tilde{\varphi} = \varphi$. Then 
$\tilde{u}(s)=u(\tilde{\varphi}(s))$ satisfies $\nabla_{J,\tilde{H}}\tilde{u}=0$ whenever $u:\IR\to M$ is a solution of 
$\nabla_{J,H}u=0$, since 
\begin{equation*}
\partial_s \tilde{u} + \nabla \tilde{H}(s,\tilde{u}) 
 = \del_s \tilde{\varphi}(s) \cdot \partial_s u + \varphi(s) \cdot \nabla H(u)\,.
\end{equation*}
For $\tilde{\eta} \in L^p(\tilde{u}^*TM)$ we find $\eta \in L^p(u^*TM)$ so that $\tilde{\eta}(s) 
= \eta(\tilde{\varphi}(s))$. Assuming that $\<F_{\tilde{u}}\tilde{\xi},\tilde{\eta}\> =0$ for all 
$\tilde{\xi}\in H^{1,p}(\tilde{u}^*TM)$, it follows that $\<F_u\xi,\eta\>=0$ for all $\xi \in H^{1,p}(u^*TM)$ 
by identifying $\tilde{\xi}(s) = \xi(\tilde{\varphi}(s))$, where $\tilde{F}_{\tilde{u}}$, $F_u$ denote the linearizations of 
$\nabla_{J,\tilde{H}}$, $\nabla_{J,H}$ at $\tilde{u},u$, respectively. The regularity of $(H,J)$ provides us with the surjectivity 
of $F_u$ at any solution $u:\IR\to M$, so that $\eta$ and therefore $\tilde{\eta}$ must vanish. $\qed$ \\

With the fixed Hamiltonian $H^{(2)}:M\to\IR$ for the cylinders we choose the Hamiltonian homotopy for the cylinders $\tilde{H}^{(2)}: \IR\times M\to\IR$ to be  
\begin{equation*} \tilde{H}^{(2)}(s,p)=\varphi(s)\cdot H^{(2)}(p), \end{equation*} 
so that $\tilde{H}^{(2)}(s,\cdot)=H^{(2)}/2$ for $s\leq -1$ and $\tilde{H}^{(2)}(s,\cdot)=H^{(2)}$. After possibly rescaling $H^{(2)}$, we can and will assume 
that both lemma 2.5 and lemma 5.1 hold with $\tau=1$ for the fixed $J$ and the chosen $H^{(2)}$, $\tilde{H}^{(2)}$, respectively. \\ 

Before we prove transversality in the next subsection, let us state the following 
analogue of theorem 2.6. Denote by $\Ju^{\tilde{H}}_N$ the domain-dependent almost complex structure on $\RS\times M$ induced by $\tilde{H}_N$. \\

\begin{theorem} Depending on the number of punctures $n$ we have the following result about the moduli spaces of $\Ju^{\tilde{H}}_N$-holomorphic 
curves in $\RS\times M$:
$ $\\
\begin{itemize}
\item $n=0$: All holomorphic spheres are constant.
\item $n=1$: Holomorphic planes do not exist.
\item $n=2$: For $T\leq 2^N$ the automorphism group $\Aut(\CP)$ acts on the parametrized moduli space $\IM^0(S^1\times M,(x^+,T),(x^-,T),\Ju^{\tilde{H}}_N)$ of holomorphic cylinders with constant finite isotropy group $\IZ_T$ and the quotient can be naturally identified with the space of gradient flow lines of $H^{(2)}$ with respect to the metric 
      $\omega(\cdot,J\cdot)$ on $M$ between the critical points $x^+$ and $x^-$ of $H^{(2)}$. In particular, we have 
      \begin{equation*}
       \sharp \IM(\RS\times M;(x^+,T),(x^-,T);\Ju^{\tilde{H}}_N) = \delta_{x^-,x^+}\,
      \end{equation*}
      since the zero-dimensional components are empty for $x^+\neq x^-$ and just contain the constant 
      path for $x^+=x^-$.
\item $n \geq 3$: For $P^+\subset P(H^{(2)}/2^N,\leq 2^N)$ and $P^-\subset P(H^{(2)}/2^{N+1},\leq 2^N)$ the action of $\Aut(\CP)$ 
      on the parametrized moduli space is free. There remains a free $S^1$-action on the moduli space, where the 
      quotient is given by
      \begin{equation*} \{(s_0,u,\uz): s_0\in\IR, u: \CP-\{\uz\} \to M: (*1),(*2)\}/\Aut(\CP)
      \end{equation*}
      with
      \begin{eqnarray*}
      &(*1):& du + X^{\tilde{H}_N}_{\uz}(z,h^0_1+s_0,u)\otimes dh_2^0 \\&& + J(u) 
              \cdot (du + X^{\tilde{H}_N}_{\uz}(z,h^0_1+s_0,u)\otimes dh_2^0)\cdot i = 0\,,\\    
      &(*2):& u\circ\psi^{\pm}_k(s,t)\stackrel{s\to\pm\infty}{\longrightarrow} x^{\pm}_k.
      \end{eqnarray*}
\end{itemize} \end{theorem}
\noindent
{\it Proof:} The proof is completely analogous to the one of theorem 2.6. Note that it follows by lemma 2.3 
that $h: \CP-\{\uz\}\to\RS$ can be identified with $(s_0,t_0)\in\RS$ and that the map $u$ now satisfies an $s_0$-dependent 
perturbed Cauchy-Riemann equation. For $n=2$ observe that by lemma 4.1 we can identify $\IM(S^1\times M;(x^+,T),(x^-,T);\Ju^{\tilde{H}}_N)$ 
with the space of all $u:\IR\to M$ satisfying  $\nabla_{J,\tilde{H}^{(2)}}u=0$, $u(s,t)\to x^{\pm}$, which following the proof of lemma 4.1 
can be identified with the space of $\tilde{u}(s)=u(\tilde{\varphi}(s))$ satisfying $\nabla_{J,H^{(2)}}u=0$. $\qed$ \\

\subsection{Transversality}

For the remaining part of this section we discuss transversality, where we again restrict ourselves to the case $N=0$. \\
\\
Since $\CR_{\Ju^{\tilde{H}}}(h,u) = (\CR h,\CR_{J,\tilde{H},s_0}u)$ with 
\begin{eqnarray*}
\CR_{J,\tilde{H},s_0}u &=& du + X^{\tilde{H}}(j,z,h^0_1+s_0,u)\otimes dh_2^0 \\&+& J(u) 
               \cdot (du + X^{\tilde{H}}(j,z,h^0_1+s_0,u)\otimes dh_2^0)\cdot i,
\end{eqnarray*}
where $X^{\tilde{H}}(j,z,s,u)$ denotes the symplectic gradient of $\tilde{H}(j,z,s,\cdot): M\to\IR$, 
it follows that the linearization $D_{h,u}$ of $\CR_{\Ju^{\tilde{H}}}$ 
is again of diagonal form. \\

It follows that for $n=2$ we get transversality from lemma 4.2 and lemma 5.1 by the special choice of $\tilde{H}^{(2)}$. \\
\\
For $n\geq 3$ let us describe the setup for the underlying universal Fredholm problem. \\
\\
As before the Cauchy-Riemann operator extends to a $C^{\infty}$-section in a Banach space bundle 
$\tilde{\EE}^{p,d}\to\BB^{p,d}\times\tilde{\H}^{\ell}$. Here $\BB^{p,d} = 
\BB^{p,d}(\RS\times M;P^+,P^-)$ denotes the manifold of maps from section 5, which is given by the product 
\begin{equation*}
   \BB^{p,d}(\RS\times M;(x_k^{\pm},T_k^{\pm})) = H^{1,p,d}_{\cst}(\Si,\IC)\times \BB^p(M;(x_k^{\pm}))\times \IM_{0,n}\,,
\end{equation*}
while the set of coherent Hamiltonian perturbations $\H^{\ell}_n(M;(H^{(k)})_{k=2}^{n-1})$ is now replaced by the 
set of coherent Hamiltonian homotopies 
\begin{equation*}
\tilde{\H}^{\ell}\;=\; \tilde{\H}^{\ell}_n(M;H;(\tilde{H}^{(k)})_{k=2}^{n-1})
\end{equation*} 
for fixed coherent Hamiltonian $H: \coprod_n \IM_{n+1}\times M\to\IR$ and $\tilde{H}^{(2)},...,\tilde{H}^{(n-1)}$. \\
Any $\tilde{H}^{(n)}\in\tilde{\H}^{\ell}$ is a $C^{\ell}$-map 
\begin{equation*}
\tilde{H}^{(n)}: \IM_{0,n+1}\times\IR\times M\to \IR,
\end{equation*} 
which extends to a $C^{\ell}$-map on $\DM_{0,n+1}\times\IR \times M$, so that 
\begin{itemize}
\item on $\bigl((\DM_{0,n+1}-\IM_{0,n+1}) \cup (\IM_{0,n}\times N_0)\bigr) \times\IR\times M$ it is given by\\ 
      $\tilde{H}^{(2)},...,\tilde{H}^{(n-1)}$, 
\item $\tilde{H}^{(n)}=H^{(n)}/2$ on $\IM_{0,n+1}\times (-\infty,-2^N)\times M$,
\item and $\tilde{H}^{(n)}=H^{(n)}$ on $\IM_{0,n+1}\times (+2^N,+\infty)\times M$,
\end{itemize}
where $N_0\subset \Si$ again denotes the fixed neighborhood of the punctures. 
It follows that the tangent space at $\tilde{H}=\tilde{H}^{(n)}\in\tilde{\H}^{\ell}$ is given by  
\begin{eqnarray*}
  T_{\tilde{H}}\tilde{\H}^{\ell}_n \;=\; \tilde{\H}^{\ell}_n(M;0;(0)_{k=2}^{n-1}).
\end{eqnarray*}

Since the linearization of $\CR_{\Ju^{\tilde{H}}}$ at $(\hb,u,j,\tilde{H})\in\BB^{p,d}\times\tilde{\H}^{\ell}$ is again 
of diagonal form,  
\begin{eqnarray*}
  &&D_{\hb,u,j,\tilde{H}} \,=\, D_j+\diag(\CR,D_{u,\tilde{H}}). \\
  && T_j\IM_{0,n} \oplus H^{1,p,d}_{\cst}(\Si,\IR^2)\oplus H^{1,p}(u^*TM)\oplus T_{\tilde{H}}\tilde{\H}^{\ell} \\
  && \to L^{p,d}(T^*\Si\otimes_{j,i}\IR^2) \oplus L^p(\Lambda^{0,1}\otimes_{j,J} u^*TM)
\end{eqnarray*}
it remains by lemma 4.2 to prove surjectivity of $D_{u,\tilde{H}}$, which is the linearization of 
the perturbed Cauchy-Riemann operator $\CR_{J,s_0}(u,\tilde{H}) = \CR_{J,\tilde{H},s_0}(u)$.  
Since the proof is in the central arguments completely similar to lemma 4.3, we just sketch the main points. \\

Assume for some $\eta\in L^q(\Lambda^{0,1}\otimes_{j,J} u^*TM)$ that $\<D_{u,\tilde{H}}(\xi,\tilde{G}),\eta\>=0$ 
for all $(\xi,\tilde{G})\in H^{1,p}(u^*TM)\oplus T_{\tilde{H}}\tilde{\H}^{\ell}$, where again $1/p+1/q=1$. From $\<\eta,D_u\xi\>=0$ for all 
$\xi$ we already know that it suffices to show that $\eta$ vanishes on an open and dense subset. 

Now observe that it follows from the same arguments used to prove lemma 4.3 that 
\begin{eqnarray*} 
 0 \,=\, \<D_{\tilde{H}} \tilde{G},\eta\> &=& 
 \int_{\Si-B} \|dh_1^0(z)\|^2\,\cdot d\tilde{G}(j,z,h^1_0(z)+s_0,u(z))\cdot\eta_1(z)\;dz 
\end{eqnarray*}
for all $\tilde{G}\in T_{\tilde{H}}\tilde{\H}^{\ell}$, where $B$ is the set of branch points of $h^0: \Si\to\RS$, we again write $\eta(z) = \eta_1(z) \otimes dh^0_1 + \eta_2(z) \otimes dh^0_2$ with 
$\eta_2(z) + J(u)\eta_1(z) = 0$ for $z\in\Si-B$ and where $d\tilde{G}(j,z,h^1_0(z)+s_0,\cdot)$ denotes the differential of 
$\tilde{G}(j,z,h^1_0(z)+s_0,\cdot): M\to\IR$. But with this we can prove as before that $\eta$ vanishes 
identically on the open and dense subset $\Si-B$. \\

Assume to the contrary that $\eta(z_0)\neq 0$, i.e., $\eta_1(z_0)\neq 0$ for some $z_0\in \Si-B$.  
As in the proof of lemma 4.3 we find $G_0 \in C^{\infty}(M)$ so that 
\begin{equation*} dG_0(u(z_0)) \cdot \eta_1(z_0) > 0. \end{equation*}

Setting $j_0:=j$, we organize all fixed maps $h_0:\Si\to\RS$ for different $j$ on $\Si$ into a map 
$h_0: \IM_{0,n+1}\to\RS$. Let $\tilde{\varphi}\in C^{\infty}(\DM_{0,n+1}\times\IR,[0,1])$ be a smooth cut-off function around 
$(j_0,z_0,h^1_0(j_0,z_0)+s_0)\in \IM_{0,n+1}\times\IR$ with $\varphi(j_0,z_0,h^1_0(j_0,z_0)+s_0)=1$ and 
$\varphi(j,z,h^1_0(j,z)+s)=0$ for $(j,z,s) \not\in U(j_0,z_0,s_0)$. Here the neighborhood 
$U(j_0,z_0,s_0) \subset \DM_{0,n+1}\times\IR$ is chosen so small that
\begin{eqnarray*} 
 U(j_0,z_0,s_0)\;\cap\; 
 \Bigl(\bigl((\DM_{0,n+1}-\IM_{0,n+1})\,\cup\,(\IM_{0,n+1}\times N_0) \bigr)\,\times\, \IR \Bigr) &=& \emptyset, \\
 U(j_0,z_0,s_0)\;\cap\; \Bigl(\DM_{0,n+1} \,\times\,\bigl((-\infty,-1)\cup(+1,+\infty)\bigr)\Bigr) &=& \emptyset,
\end{eqnarray*} 
and $dG_0(z,u(z))\cdot\eta_1(z) \geq 0$ for all $(z,j,h^1_0(j,z)+s) \in U(j_0,z_0,s_0)$. \\

Defining $\tilde{G}: \DM_{0,n+1}\times\IR\times M \to \IR$ by $\tilde{G}(j,z,s,p):=\varphi(j,z,s) \cdot G_0(p)$, 
this leads to the desired contradiction since we found $\tilde{G} \in T_{\tilde{H}} \tilde{\H}^{\ell}=
\tilde{\H}^{\ell}_n(M;0;0,...,0)$ with 
\begin{eqnarray*}
\<D_{\tilde{H}}\cdot \tilde{G},\eta\> = 
\int_{\Si-B} \|dh_1^0(z)\|^2\,\cdot\,d\tilde{G}(j_0,z,h^1_0(j_0,z)+s_0,u(z))\cdot \eta_1(z) \;dz > 0. 
\end{eqnarray*} 

We have shown that the universal moduli space 
$\IM(\RS\times M;P^+,P^-;\Ju^H;(\Ju^{\tilde{H},(k)})_{k=2}^{n-1})$ is again transversally cut out by the Cauchy-Riemann operator 
$\CR_J$. Further it follows by the same arguments as in section 4 that we can choose a (smooth) coherent Hamiltonian 
homotopy $\tilde{H}: \coprod_n \IM_{0,n+1} \times \IR\to C^{\infty}(M)$ such that for all $N\in\IN$ and $P^+,P^-$ 
the moduli spaces $\IM(\RS\times M;P^+,P^-;\Ju^{\tilde{H}}_N)$ are transversally cut out by the Cauchy-Riemann operator. 
  
\section{Contact homology}

\subsection{Chain complex} 

The contact homology of $S^1\times M$ equipped with the stable Hamiltonian structure $(\omega^H,\lambda^H)$ is defined as the homology of a 
differential graded algebra $(\IA,\del)$, which is generated by closed orbits of the Reeb vector field $R^H$ and whose differential counts 
$\Ju^H$-holomorphic curves with one positive puncture. As in [EGH] we start with assigning to any 
$(x,T) \in P(H)$, which is {\it good} in the sense of [BM], a graded variable $q_{(x,T)}$ with 
\begin{equation*} \deg q_{(x,T)} = \dim M/2 -2 + \mu_{CZ}(x,T). 
\end{equation*}
Here $\mu_{CZ}$ denotes the Conley-Zehnder index for $(x,T)$, which is defined as in [EGH] 
after fixing a basis for $H_1(S^1\times M)$ and choosing a spanning surface between the orbit $(x,T)$
and suitable linear combinations of these basis elements. Note that in the corresponding definition in [EGH] one adds $m-3$, 
where $m$ denotes the complex dimension of $\RS\times M$. Further we assume, as in [EGH], that $H_1(S^1\times M)$ and hence $H_1(M)$ is torsion-free, where we use that the torsion-freedom of $H_*(S^1)$ also yields the Kuenneth formula for $H_*(S^1\times M)$. Let 
\begin{equation*} \IQ[H_2(S^1\times M)] = \{\sum q(A) e^A: A\in H_2(S^1\times M), q(A) \in \IQ\} \end{equation*} 
be the group algebra generated by $H_2(S^1\times M) \cong H_2(M) \oplus (H_1(S^1)\otimes H_1(M))$ with grading given by
\begin{equation*} \deg e^A = \<c_1(TM),A\>. \end{equation*}
Since $c_1(TM)$ clearly vanishes on $H_1(S^1)\otimes H_1(M)$ we can and will work with the reduced group ring $\IQ[H_2(M)]$. With this let $\IA_*$ be the graded commutative algebra of polynomials in the formal variables $q_{(x,T)}$ assigned to good periodic orbits with coefficients in $\IQ[H_2(M)]$. Let $C_*$ be the vector space over $\IQ$ freely generated by the graded variables $q_{(x,T)}$, which naturally splits, $C_* = \bigoplus_T C_*^T$ with 
$C_*^T$ generated by the good orbits of integer period $T$. Since $C_*$ is graded, we can define a graded symmetric algebra 
$\SS(C_*)$ and it follows that 
\begin{equation*} \IA_* = \SS(C_*)\otimes \IQ[H_2(M)]. \end{equation*}

For the following we assume that all occuring periodic orbits are good. Note that to any holomorphic curve in $\IM(S^1\times M;P^+,P^-;\Ju^H)$ we assign as in [EGH] a homology class $A\in H_2(S^1\times M)$ after fixing a basis for $H_1(S^1\times M)$ and choosing spanning surfaces between the asymptotic orbits in $P^+,P^-\subset P(H)$ and suitable linear combinations of these basis elements. Requiring that the differential $\del: \IA \to \IA$ satisfies a graded Leibniz rule, it is defined by (see [EGH],p.621)
\begin{equation*}
   \del q_{(x_0,T_0)} = \sum_{P^-,A} \sharp \IM_A(S^1\times M;P^+,P^-;\Ju^H)/\IR\;\,q_{(x_1^-,T_1^-)} ... q_{(x_n^-,T_n^-)}\;e^A, 
\end{equation*}
where $\IM_A(S^1\times M;P^+,P^-;\Ju^H)$ denotes the one-dimensional component of the moduli space of holomorphic curves with $P^+ = \{(x_0,T_0)\}$ but arbitary orbit set $P^- = \{(x_1^-,T_1^-),...,(x_n^-,T_n^-)\}$ representing the class $A\in H_2(M) \cong H_2(S^1\times M)/(H_1(S^1)\otimes H_1(M))$. \\

For $(T_1,...,T_n)\in\IN^n$ let $\IA^{(T_1,...,T_n)}$ denote the subspace of $\IA$ spanned by monomials 
$q_{(x_1,T_1)}\, ... \,q_{(x_n,T_n)}$, 
\begin{equation*} 
\IA^{(T_1,...,T_n)} = \SS^{(T_1,...,T_n)}(C_*) \otimes \IQ[H_2(M)] 
\end{equation*}
with
\begin{equation*} 
\SS^{(T_1,...,T_n)}(C_*) = \SS(C_*^{T_1} \otimes ... \otimes C_*^{T_n}),
\end{equation*}
where $\SS$ denotes the projection from the tensor to the symmetric algebra, in particular, $\IA^{(T_1,...,T_n)}$ does not depend on the 
ordering of the $T_1,...,T_n$. Since $\sharp \IM(S^1\times M;P^+,P^-;\Ju^H)/\IR = 0$ for $T_1^-+...+T_n^-\neq T_k$ by lemma 1.1.3, 
it follows that the differential $\del$ respects the splitting 
\begin{equation*} 
  \IA = \bigoplus_{T\in\IN} \IA^T,  
\end{equation*} 
where $\IA^T = \bigoplus_{T_1+...+T_n=T} \IA^{(T_1,...,T_n)}$. 

\subsection{Proof of the main theorem} 

In what follows we use our results about holomorphic curves in $\RS\times M$ to prove the main theorem. 
At first we compute $H_*(\IA^{\leq 2^N},\del) = \bigoplus_{T\leq 2^N} H_*(\IA^T,\del)$ using our results about moduli spaces of holomorphic curves in $\RS\times M$ in theorem 2.6 together with the transversality results. \\

With the fixed almost complex structure $J$ on $M$ let $H: \coprod \IM_{0,n+1} \to C^{\infty}(M)$ be a coherent Hamiltonian 
perturbation as before, in particular, $H^{(2)}$ satisfies lemma 2.5 with $\tau=1$. Following corollary 4.4 we further assume 
that $H$ is chosen such that transversality holds for all moduli spaces $\IM(S^1\times M;P^+,P^-;\Ju^{H/2^N})$, 
$P^{\pm}\subset P(H^{(2)}/2^N,\leq 2^N)$, simultaneously for all $N\in\IN$. Together with theorem 2.6 it then follows that 
for defining the algebraic invariants we only have to count gradient flow lines of the function $H^{(2)}$ on $M$ with respect 
to the metric $g_J=\omega(\cdot,J\cdot)$ on $M$. For $N \in \IN$ let $(\IA_{N},\del_{N})$ denote the differential algebra for the domain-dependent Hamiltonian $H/2^N: \coprod \IM_{0,n+1} \to C^{\infty}(M)$ and the fixed almost complex structure $J$ on $M$. 
For the computation of the contact homology subcomplex we use special choices for the basis elements in $H_1(S^1\times M)$ 
and the spanning surfaces as follows: Choose a basis for $H_1(S^1\times M)=H_1(S^1)\oplus H_1(M)$ containing the 
canonical basis element $[S^1]$ of $H_1(S^1)$, which is represented by the circle $(x^*,1): S^1\to S^1\times M$, 
$t\mapsto(t,x^*)$ for some point $x^*\in M$. For any periodic orbit $(x,T)\in P(H^{(2)}/2^N,\leq 2^N)$ we have 
$[(x,T)]=T[S^1]\in H_1(S^1\times M)$, since $x$ is a constant orbit in $M$, and we naturally specify 
a spanning surface $S_{(x,T)}$ between $(x,T)$ and the $T$-fold cover of $(x^*,1)$ by choosing a path 
$\gamma_x: [0,1]\to M$ from $x^*$ to $x$ and setting $S_{(x,T)}: S^1\times [0,1] \to S^1\times M$, 
$S_{(x,T)}(t,r)=(Tt,\gamma_x(r))$. \\
\\
\begin{lemma} Let $HM_*=HM_*(M,-H^{(2)},g_J;\IQ)$ denote the 
Morse homology for the Morse function $-H^{(2)}$ and the metric $g_J=\omega(\cdot,J\cdot)$ on $M$ with 
rational coefficients. Then we have
\begin{equation*}
  H_*(\IA^{\leq 2^N}_{N},\del_{N}) = \SS^{\leq 2^N}(\bigoplus_{\IN}HM_{*-2})\otimes \IQ[H_2(M)]\,,
\end{equation*}
where
\begin{equation*}
\SS^{\leq 2^N}(\bigoplus_{\IN}HM_{*-2})= \bigoplus_{T_1+...+T_n\leq 2^N} \SS^{(T_1,...,T_n)}(\bigoplus_{\IN} HM_{*-2}).
\end{equation*} \end{lemma}
\noindent
{\it Proof:} For the grading of the $q$-variables we have  
\begin{equation*}
\deg q_{(x,T)} = \dim M/2 - 2 + \mu_{CZ}(x,T) = \ind_{-H}(x) - 2, 
\end{equation*}
when we choose a canonical trivialization of $TM$ over $(x^*,1)$ and extend it over the spanning surfaces 
to a canonical trivialization over $(x,T)$, i.e., the map $\Theta: S^1 \times \IR^{2m} \to x^*TM = S^1 \times T_x M$ is 
independent of $S^1$. It follows that $C^T_*$ agrees with the chain group $CM_{*-2}$ for the Morse homology  
for $T\leq 2^N$ and therefore 
\begin{equation*}
\IA^{\leq 2^N}_{N} = \SS^{\leq 2^N}(\bigoplus_{\IN} CM_{*-2})\otimes \IQ[H_2(M)]\,.
\end{equation*}
Here it is important to observe that any $(x,T)\in P(H^{(2)}/2^N,\leq 2^N)$ is indeed
good in the sense of [BM]: note that it follows from $\mu_{CZ}(x,T) = \ind_{-H}(x) - \dim M/2$  
that $\mu_{CZ}(x,T)$ has the same parity for all $T\leq 2^N$. \\

On the other hand it follows from theorem 2.6 that the differential $\del:\IA\to\IA$ indeed agrees with the differential in Morse homology. Further it follows from the above choice of spanning surfaces that they all represent the trivial class $A \in H_2(M)=H_2(S^1\times M)/(H_1(S^1)\otimes H_1(M))$: Indeed, letting $u$ denote the gradient flow line between 
$x_0$ and $x$ it follows that $u$ represents the class $A = T[S^1] \otimes [\gamma_{x_0}\sharp u\sharp -\gamma_x] 
\in H_1(S^1)\otimes H_1(M)$. Using the theorem of Kuenneth we hence in fact have 
\begin{eqnarray*}
 H_*(\IA_N^{(T_1,...,T_n)},\del) 
 &=& H_*(\SS^{(T_1,...,T_n)}(\bigoplus_{\IN}CM_{*-2}) \otimes \IQ[H_2(M)], \del) \\
 &=& \SS^{(T_1,...,T_n)}(H_*(\bigoplus_{\IN}CM_{*-2},\del^{\textrm{Morse}}))\bigr) \otimes \IQ[H_2(M)]\\
 &=& \SS^{(T_1,...,T_n)} \bigl(\bigoplus_{\IN} HM_{*-2}\bigr) \otimes \IQ[H_2(M)]
\end{eqnarray*}
and the claim follows. $\qed$ \\

With this we can now complete the proof of the main theorem by using theorem 5.2 and the transversality result of section five. \\

To this end choose a coherent Hamiltonian homotopy $\tilde{H}: \coprod_n \IM_{0,n+1} \times \IR\to C^{\infty}(M)$ 
as in section five, i.e., with $\tilde{H}(j,z,s,p)=H(j,z,p)/2$ for small $s$ and $\tilde{H}(j,z,s,p)=H(j,z,p)$ for 
large $s$ such that for all $N\in\IN$ and $P^+,P^-$ the moduli spaces $\IM(\RS\times M;P^+,P^-;\Ju^{\tilde{H}}_N)$ 
are transversally cut out. Let $\Ju^{\tilde{H}}_N$ denotes the coherent non-cylindrical almost complex structure on 
$\RS\times M$ induced by $J$ and $\tilde{H}/2^N$. \\

Let $\Psi_N: (\IA_N,\del_N) \to (\IA_{N+1},\del_{N+1})$ be the chain homotopy, defined as in [EGH], by counting holomorphic curves with one positive puncture and an arbitrary number of negative punctures in the resulting almost complex manifold $(\RS\times M,\Ju^{\tilde{H}}_N)$ with cylindrical ends. Then it follows from theorem 5.2 that the restriction $\Psi^T_N: (\IA^T_N,\del_N) \to (\IA^T_{N+1},\del_{N+1})$ is the identity for $T\leq 2^N$, since again all curves with three or more punctures 
come in $S^1$-families and all zero-dimensional cylindrical moduli spaces just consist of trivial gradient flow lines. \\


\begin{thebibliography}{10000000}

\bibitem[BC]{BC} Bourgeois, F. and  Colin, V., Homologie de contact des varieties toroidales. {\it Geom. Topol.} {\bf 9}(2005), 299-313.

\bibitem[BEHWZ]{BEHWZ} Bourgeois, F., Eliashberg, Y., Hofer, H., Wysocki, K. and Zehnder, E., Compactness results in 
      symplectic field theory. {\it Geom. Topol.} {\bf 7}(2003), 799-888. 

\bibitem[BM]{BM} Bourgeois, F. and  Mohnke, K., Coherent orientations in symplectic field theory. {\it Math. Z.} {\bf 248}(2003), 123-146.

\bibitem[CC]{CC} Cotton-Clay, A., Symplectic Floer homology of area-preserving surface symplectomorphisms. 
      ArXiv preprint (0807.2488), 2008.

\bibitem[CL]{CL} Cieliebak, K. and Latschev, J.: The role of string topology in symplectic field theory. ArXiv preprint (0706.3284), 2007.

\bibitem[CM1]{CM1} Cieliebak, K. and Mohnke, K., Symplectic hypersurfaces and transversality for Gromov-Witten theory. {\it J. Symp. Geom.} {\bf 5}(2007), 281-356.
      .
\bibitem[CM2]{CM2} Cieliebak, K. and Mohnke, K., Compactness of punctured holomorphic curves. {\it J. Symp. Geom.} {\bf 3}(2005), 589-654.
 
\bibitem[EGH]{EGH} Eliashberg, Y., Givental, A. and Hofer, H., Introduction to symplectic field theory. {\it 
      Geom. Funct. Anal. 2000 Visions in Mathematics special volume, part II} (2000), 560-673.

\bibitem[EKP]{EKP} Eliashberg, Y., Kim, S. and Polterovich, L., Geometry of contact transformations and domains: 
      orderability vs. squeezing. {\it Geom. Topol.} {\bf 10}(2006), 1635-1748.

\bibitem[F]{F} Fabert, O., Counting trivial curves in rational symplectic field theory. ArXiv preprint (0709.3312), 2007.

\bibitem[FO]{FO} Fukaya, K. and Ono, K., Arnold conjecture and Gromov-Witten invariants for general symplectic manifolds. {\it 
      Topology} {\bf 38}(1999), 933-1048.

\bibitem[H]{H} Hofer, H., Polyfolds and a general Fredholm theory. ArXiv preprint (0809.3753), 2008. 

\bibitem[K]{K} van Koert, O., Contact homology of Brieskorn manifolds. {\it Forum Math.} {\bf 20}(2008), 317-339.

\bibitem[LT]{LT} Li, J. and Tian, G., Virtual moduli cycles and Gromov-Witten invariants of general symplectic manifolds. {\it 
      First Int. Press Lect. Ser. I} (1998), 47-83.
 
\bibitem[LiuT]{LiuT} Liu, G. and Tian, G., Floer homology and Arnold conjecture. {\it J. Diff. Geom.} {\bf 49}(1998), 1-74.

\bibitem[MDSa]{MDSa} McDuff, D. and Salamon, D.A., {\it $J$-holomorphic curves and symplectic topology.} AMS Colloquium Publications, 
      Providence RI, 2004.

\bibitem[Sch]{Sch} Schwarz, M., {\it Cohomology operations from $S^1$-cobordisms in Floer homology.} Ph.D. thesis, 
      Swiss Federal Inst. of Techn. Zurich, Diss. ETH No. 11182, 1995.

\bibitem[SZ]{SZ} Salamon, D.A. and Zehnder, E., Morse theory for periodic solutions of Hamiltonian systems and 
      the Maslov index. {\it Comm. Pure Appl. Math.} {\bf 45}(1992), 1303-1360.

\bibitem[Y1]{Y1} Yau, M.L., Cylindrical contact homology of subcritical Stein fillable contact manifolds. {\it Geom. and Topol.} {\bf 8}(2004), 
        1243-1280. 

\bibitem[Y2]{Y2} Yau, M.L., Vanishing of the contact homology of overtwisted contact 3-manifolds. {\it Bull. Inst. Math. Acad. Sinica} {\bf 1}(2006), 211-229.

   
\end{thebibliography}
\end{document}